\documentclass[11pt]{amsart}

\usepackage{amsmath, amsthm, amssymb, amsfonts, mathtools, mathrsfs, tikz-cd}
\usepackage{hyperref}
\usepackage[backend = biber, style = alphabetic]{biblatex}

\theoremstyle{plain}
\newtheorem{theorem}{Theorem}[section]
\newtheorem{lemma}[theorem]{Lemma}
\newtheorem{corollary}[theorem]{Corollary}
\newtheorem{proposition}[theorem]{Proposition}

\theoremstyle{definition}
\newtheorem{definition}[theorem]{Definition}

\theoremstyle{remark}
\newtheorem{remark}[theorem]{Remark}

\numberwithin{equation}{section}

\newcommand{\C}{\mathbb{C}}
\newcommand{\R}{\mathbb{R}}
\newcommand{\Q}{\mathbb{Q}}

\newcommand{\Z}{\mathbb{Z}}

\renewcommand{\P}{\mathbb{P}}
\renewcommand{\O}{\mathcal{O}}

\renewcommand{\H}{\mathcal{H}}

\renewcommand{\d}{\mathrm{d}}
\renewcommand{\i}{\mathrm{i}}

\newcommand{\tcr}[1][f]{\bar{\partial}_{#1}}
\newcommand{\tcrc}[1][f]{\partial_{#1}}
\newcommand{\tcrd}[1][f]{\bar{\partial}_{#1}^\dagger}
\newcommand{\tcrcd}[1][f]{\partial_{#1}^\dagger}

\newcommand{\ctr}[1]{\iota_{\partial_{#1}}}

\title{Hodge Theory of Novikov Type Landau--Ginzburg Models}

\author{Huijun Fan}
\address{School of Mathematics and Statistics, Wuhan University, Wuhan, Hubei, China}
\email{fanhuijun@whu.edu.cn, fanhj@math.pku.edu.cn}
\thanks{supported by ZYYJQT202600011, National Key R\&D Program of China (Grant No. 2020YFE0204200), NSFC(11831017,11890660, 11890661) and NSFC-RSF 1201101428, NSFC-RFBR 11661131005.}
\author{Jieran Wang}
\address{School of Mathematics, Sun Yat-sen University, Guangzhou, Guangdong, China}
\email{wangjr88@mail.sysu.edu.cn}

\begin{document}

\begin{abstract}
We establish a strong Hodge decomposition theorem for Novikov type Landau--Ginzburg models, that is, triples $(M,g,\theta)$ consisting of a non-compact complete K\"ahler manifold $(M,g)$ and a non-trivial closed holomorphic one-form $\theta$ on $M$.
Under the well-tameness condition and under the assumption that $M$ is Stein, we prove that the space of $L^{2}$-forms on $M$ admits an orthogonal decomposition into the space of harmonic forms of the twisted Laplacian $\Delta_{\theta}$ and the images of the twisted Dolbeault operators $\bar\partial+\theta\wedge$ and $\bar\partial^{\dagger}+\iota_{\bar\theta^{\sharp}}$; moreover, the harmonic space is concentrated in degree $n=\dim_{\C}M$ and is isomorphic to the Milnor rings of $\theta$.

As applications, we endow the complement of a hyperplane arrangement with a complete cylindrical K\"ahler metric, and we show that the differential of a master function with nonresonant weights is well-tamed with respect to this metric. After passing to a log resolution in the non-normal-crossing case; we then apply the results to the Yang--Yang functional of the Gaudin model, whose critical points are the Bethe roots: the top-degree harmonic space is described by the corresponding Milnor rings, and the number of Bethe roots is given by the Varchenko--Orlik--Terao lattice formula.
\end{abstract}

\maketitle

\section{Introduction}\label{sec:introduction}

Landau--Ginzburg (LG) models have been studied for a long time by physicists as fundamental objects of topological field theories and mirror symmetry \parencite{CecottiVafa1991}.
They are regarded as the mirror dual of Fano manifolds: the $A$-model of a Fano manifold is mirror to the $B$-model of the corresponding LG superpotential \parencite{BarannikovKontsevich1998, HoriVafa2000}.
From the perspective of singularity theory, LG models appear naturally as the targets of FJRW theory for non-degenerate quasi-homogeneous polynomials \parencite{FanJarvisRuan2013} and as the natural geometric settings for primitive forms, Frobenius manifolds and $tt^{*}$-geometry \parencite{Dubrovin1990ModernIII, Dubrovin1996}.

A systematic differential-geometric study of LG models was initiated by the first author of \parencite{fan2011}.
Despite the variety of definitions appearing in the literature, one may adopt the \emph{K\"ahlerian} version that parallels the definition of a compact K\"ahler manifold: an LG model $(M,g,f)$ is a connected, non-compact complete K\"ahler manifold together with a non-trivial holomorphic function $f:M\to\C$.
The associated \emph{twisted Dolbeault operators}
\begin{equation*}
\tcr[f]=\bar\partial+\partial f\wedge ,\qquad \tcrc[f]=\partial+\bar\partial\bar f\wedge
\end{equation*}
and their $L^{2}$-adjoints give rise to a twisted Hodge--Laplace operator $\Delta_{f}$.
When $M$ is compact, the resulting Hodge theory reproduces the classical picture.
When $M$ is non-compact, however, one must control the behaviour of $|\partial f|$ at infinity in order to obtain finiteness of the harmonic space and the strong Hodge decomposition.
Several notions of \emph{tameness} have been proposed for this purpose.
The \emph{fundamental tameness} (only finitely many critical points and full range of $f$) and \emph{strong tameness} were discussed in \parencite{fan2011}; the \emph{strong tameness} is also called the elliptic condition in the early paper \parencite{KlimekLesniewski1991}. The \emph{well-tameness} adopted below is weaker than \emph{strong tameness} and was singled out by \parencite{DaiYan2022} in their study of Witten deformation for non-compact manifolds with bounded geometry.

The present paper addresses the \emph{Novikov type} generalisation, in which the holomorphic function $f$ is replaced by a non-trivial closed holomorphic $1$-form $\theta$ on $(M,g)$.
Such triples $(M,g,\theta)$ arise naturally in several situations: as section-bundle systems \parencite{fan2011}, as the deformation parameters of complex structures appearing in Hodge-theoretic aspects of mirror symmetry \parencite{KatzarkovKontsevichPantev2008}, and as the input data for Witten deformations that are not globally exact.
Formally one defines
\begin{equation*}
\tcr[\theta]=\bar\partial+\theta\wedge ,\qquad \tcrc[\theta]=\partial+\bar\theta\wedge ,
\end{equation*}
replacing the exterior product $\partial f\wedge$ by wedging with $\theta$.
The weak Hodge decomposition (Equation \eqref{weak_Hodge} in the main text) is then automatic on complete K\"ahler manifolds, but again the strong decomposition requires a quantitative condition at infinity.
The relevant notion is the following \emph{well-tameness} condition, introduced by \parencite{DaiYan2022} and adopted as Definition~\ref{def:well-tamed} in Section~\ref{sec:main_theorem}: there exist a compact set $K\subset M$ and positive constants $C_{1},C_{2}$ such that outside $K$ one has $|\theta|>C_{1}$ and $|\theta|^{2}-C_{2}\,|\nabla\theta|>0$.

The main theorem of this paper is the following; its proof will be given in Section~\ref{sec:main_theorem}.

\begin{theorem}
Let $(M,g,\theta)$ be a well-tamed Novikov type LG model and assume that $M$ is Stein.
Then the strong Hodge decomposition
\begin{equation*}
L^{2}(\Lambda^{\bullet}(M)) = \H^{\bullet}_{\theta}\oplus \mathrm{im}(\tcr[\theta])\oplus \mathrm{im}(\tcrd[\theta])
\end{equation*}
holds, and the harmonic space is concentrated in the top degree and decomposes as
\begin{equation*}
\H^{n}_{\theta}\cong\bigoplus_{x\in \mathrm{Crit}(\theta)} \mathcal{R}_{\theta,x}.
\end{equation*}
\end{theorem}

Theorem~\ref{main_theorem_1} extends the strong Hodge theorem for function-type LG models in \parencite{FanFang2021} to closed holomorphic one-forms.
Its proof combines the homotopy argument with a reduction to the holomorphic Koszul complex.

\begin{remark}
Well-tameness is stronger than fundamental tameness. Thus the critical set is a compact analytic subset of $M$, and the Stein assumption makes it finite.
The cohomological reduction in Section~\ref{sec:main_theorem} then identifies the twisted cohomology with the local Milnor rings supported at these critical points.
\end{remark}

We apply this framework to hyperplane arrangement complements and to the Yang--Yang functional.\
 \\
\paragraph{Application~I: Hyperplane arrangement complements.}
Let $\mathcal A=\{H_1,\dots,H_m\}$ be an arrangement in $\P^n$, with $H_i=\{f_i=0\}$ and $f_i(Z)=\sum_{j=0}^n c_{ij}Z_j$, and set $M=\P^n\setminus\bigcup_iH_i$.
Section~\ref{sec:YY} recalls the intersection poset $L(\mathcal A)$ and the critical-point counting formula for the master function
\begin{equation*}
 F=\sum_i a_i\log f_i,\qquad
 \theta=\partial F=\sum_i a_i\frac{\partial f_i}{f_i},
 \qquad \sum_i a_i=0.
\end{equation*}
For nonresonant weights, the Varchenko--Orlik--Terao formula \parencite[Lemma~2.8]{Varchenko2015CriticalSet} gives
\begin{equation*}
 \sum_{p\in Z(\theta)}\mu_p=(-1)^n\chi(M)
 =(-1)^n\sum_{X\in L(\mathcal A)}\mu(X)(\dim X+1),
\end{equation*}
where $\mu_p$ is the local Milnor number and $\mu(X)$ is the M\"obius function of the intersection poset.

We construct a complete K\"ahler metric on $M$ from
\begin{equation*}
 \omega=\omega_{FS}+\frac{\i}{2}\sum_{i=1}^m
            \partial\bar\partial(\rho_i^2),\qquad
 \rho_i=C_i-\log\|f_i\|_{FS}^2\ge1,\qquad
 C_i=1+\log\left(\sum_{j=0}^n|c_{ij}|^2\right).
\end{equation*}
Its expression
\begin{equation*}
 \omega=(1+Q)\omega_{FS}+\i\sum_i\partial\rho_i\wedge\bar\partial\rho_i,
 \qquad Q=\sum_i\rho_i,
\end{equation*}
contains both normal logarithmic terms and a growing background factor $1+Q$.
For a simple normal crossing arrangement, weighted estimates in a logarithmic frame give
\begin{equation*}
 |\nabla\theta|_g\le C_0,\qquad
 |\theta|_g^2\ge c_0>0
\end{equation*}
outside a compact subset.
Choosing $C_2$ sufficiently small establishes well-tameness.

Applying Theorem~\ref{main_theorem_1} to these models, we obtain Proposition~\ref{prop:well-tamed-normal-crossing} and its corollary, where the top-degree harmonic dimension is identified with the total Milnor number. 

Section~\ref{sec5.4} discusses the passage to a log resolution for a general arrangement, where the boundary residues are the corresponding sums of hyperplane weights. The LG model of Yang-Yang functional is an example that in general the natural boundary is formed by non-normal crossing divisors. Therefore before applying Theorem~\ref{main_theorem_1}, we need to do log resolution.\
 \\
\paragraph{Application~II: The Yang--Yang functional and the Bethe ansatz.}
For a semisimple Lie algebra $\mathfrak g$ with marked points $z=(z_1,\dots,z_N)$, the Yang--Yang functional $\Phi(t;z)$ is a multivalued function whose critical-point equations in the Bethe variables $t$ are the Bethe equations of the Gaudin model \parencite{NEKRASOV201169, Liu_2021, HuLiu2015}.
Its differential $\partial\Phi$ is a closed holomorphic logarithmic one-form on the complement $M_z$ of a collision arrangement.
Section~\ref{sec:YY-functional} describes this arrangement and the boundary residues associated with collisions. The main result is Proposition~\ref{prop:YY-finiteness}. It gives the following conclusions. When $(M_z,g_z,\partial\Phi)$ satisfies the nonresonant condition, its top-degree harmonic space is described by the local Milnor rings at the Bethe roots. The total Milnor number is given by \eqref{eq:counting-lattice}; for weights with nondegenerate critical points, this also counts the distinct Bethe roots.\
 
\begin{remark}[Comparison with the circle-valued Morse theory of Kohno--Pajitnov.]
Kohno and Pajitnov \parencite{KohnoPajitnov2015} study an essential affine arrangement with complement $M=\C^n\setminus\bigcup_iH_i$ using positive real weights $\alpha_i$ and the function
\begin{equation*}
 f_\alpha=\prod_i|f_i|^{\alpha_i},\qquad
 \omega_\alpha=\sum_i\alpha_i\frac{\partial f_i}{f_i},
 \qquad \d\log f_\alpha=\operatorname{Re}\omega_\alpha.
\end{equation*}
For generic positive weights, $f_\alpha$ is Morse and all its critical points have index $n$.
For sufficiently small $\varepsilon>0$, $M$ has the homotopy type of the level manifold $\{f_\alpha=\varepsilon\}$ with $|\chi(M)|$ cells of dimension $n$ attached; for positive weights of rank one over $\Q$, this level manifold fibers over $S^1$.
They also show that the Novikov homology for any positive homomorphism $\xi:\pi_1(M)\to\R$ is concentrated in degree $n$ and is free of rank $|\chi(M)|$ over the corresponding Novikov ring.
Our approach studies the holomorphic logarithmic one-form with complex weights through the $L^2$ Hodge theory of $\tcr[\theta]$.
The invariant $|\chi(M)|$ appears in their cell and Novikov-homology counts, and as the total Milnor number in \eqref{eq:counting-lattice}.
\end{remark}

\paragraph{Organisation of the paper.}
Section~\ref{sec:LG} introduces Novikov type LG models and their twisted Dolbeault operators.
Section~\ref{sec:main_theorem} proves Theorem~\ref{main_theorem_1} through a homotopy argument and a reduction to the holomorphic Koszul complex.
Section~\ref{sec:YY} recalls the intersection poset and master-function counting formula, constructs a complete K\"ahler metric, and proves well-tameness for simple normal crossing arrangements.
Section~\ref{sec:YY-functional} studies the Yang--Yang functional, the collision arrangement and the Bethe equations.

\section{Novikov Type Landau--Ginzburg Models}\label{sec:LG}

\subsection{Definitions}

In this section we introduce Novikov type Landau--Ginzburg models and collect the analytic prerequisites for the rest of the paper.
Since the closed holomorphic one-form $\theta$ is locally exact, all local constructions of singularity theory, such as critical points, the Milnor ring and the Milnor numbers, carry over verbatim to the closed-form setting; globally, the twisting of the Dolbeault operators by $\theta$ produces a complex whose weak Hodge decomposition already holds on any complete K\"ahler manifold.
We begin with the definition.

\begin{definition}
	The tuple $(M,g,\theta)$ is called a Novikov type LG model if
	\begin{enumerate}
		\item $(M,g)$ is a non-compact complete K\"ahler manifold.
		\item $\theta$ is a nontrivial closed holomorphic $1$-form on $M$.
	\end{enumerate}
\end{definition}

\begin{remark}
	The Novikov type LG model is also a special case of section-bundle systems (see \parencite{fan2011}). 
	We note that there always exists a covering $p: \tilde{M} \to M$ such that $p^*\theta$ is exact. 
	Therefore, a Novikov type LG model can be viewed as an LG model equipped with a free action of a discrete group.
\end{remark}

Let $(M,g,\theta)$ be a Novikov type LG model and let $M$ be of dimension $n$. 
Since $\theta$ is locally exact, we can import the constructions from singularity theory. 
Let $f$ be a local primitive of $\theta$. 
We call a point $p \in M$ a (Morse) critical point of $\theta$ if $p$ is a (Morse) critical point of $f$. 
We denote the set of critical points of $\theta$ by $\mathrm{Crit}(\theta)$. 
We denote the sheaf of Jacobian ideal by
\begin{equation}
	\mathcal{J}_\theta = \left(\theta_1,\ldots, \theta_n\right),
\end{equation}
where $\theta_i$ is the component of $\theta$ under some local coordinates.
We define the sheaf of Milnor ring of $\theta$ by
\begin{equation}
	\mathcal{R}_\theta = \O_M/\mathcal{J}_\theta.
\end{equation}
Let $p$ be a critical point of $\theta$. 
The stalk of $\mathcal{R}_\theta$ at $p$ is called the Milnor ring of $\theta$ at $p$ and is denoted by $\mathcal{R}_{\theta,p}$.
The dimension of $\mathcal{R}_{\theta,p}$ as a $\C$-vector space is called the Milnor number of $\theta$ at $p$ and is denoted by $\mu_p$.
We can also define the global Milnor number of $\theta$ by
\begin{equation}
    \mu_\theta = \begin{cases}
        \sum_{p \in \mathrm{Crit}(\theta)} \mu_p. & \text{if $\mathrm{Crit}(\theta)$ is finite}, \\
        \infty. & \text{otherwise}.
    \end{cases}
\end{equation}

\subsection{Twisted Dolbeault operators}

We define the twisted $\bar\partial$ operators by
\begin{align}
	 & \tcr[\theta] = \bar\partial + \theta\wedge   \\
	 & \tcrc[\theta] = \partial + \bar\theta\wedge.
\end{align}
We denote by $\tcrd[\theta]$ and $\tcrcd[\theta]$ the adjoints of $\tcr[\theta]$ and $\tcrc[\theta]$ with respect to the K\"ahler metric,
\begin{align}
	 & \tcrd[\theta] = \bar\partial^\dagger + \iota_{\bar\theta^\sharp} \\
	 & \tcrcd[\theta] = \partial^\dagger + \iota_{\theta^\sharp}.
\end{align}
We define the twisted Laplacian by
\begin{equation}
	\Delta_\theta = (\tcr[\theta] + \tcrd[\theta])^2.
\end{equation}

Since $\theta$ is locally exact, these operators share the same local properties as operators from classical LG models. 
These properties include commutative relations, K\"ahler-Hodege identities and relations with Hodge $*$ operators.
Precise description can be found in \parencite{fan2011}.

We can denote $\tcr[\theta]$, $\tcrc[\theta]$, $\tcrd[\theta]$, $\tcrcd[\theta]$ as their own closure on $L^2(\Lambda^\bullet (M))$, and denote $\Delta_\theta$ as its own Friedrichs extension on $L^2(\Lambda^\bullet (M))$.
The following property is obvious:
\begin{proposition}
	Let $\psi \in \mathrm{Dom}(\Delta_\theta)$.
	Then $\Delta_\theta \psi = 0$ if and only if $\tcr[\theta] \psi = 0$ and $\tcrd[\theta] \psi = 0$.
\end{proposition}
Now we can define the space of $\Delta_\theta$-harmonic forms by
\begin{equation}
	\mathcal{H}^\bullet_\theta = \ker \Delta_\theta.
\end{equation}
The proposition above implies the following weak Hodge decomposition
\begin{equation}\label{weak_Hodge}
	L^2(\Lambda^\bullet(M)) = \H_\theta^\bullet \oplus \overline{\mathrm{im} (\tcr[\theta])} \oplus \overline{\mathrm{im} (\tcrd[\theta])}.
\end{equation}

\section{Main Theorem and Proof}\label{sec:main_theorem}
\subsection{Main Theorem}

This section is devoted to the statement and the proof of the main theorem.
The proof proceeds in three steps.
First, a homotopy argument identifies the $L^{2}$-cohomology of the $\tcr[\theta]$-complex with its ordinary and compactly supported cohomologies.
Second, a $\tcr[\theta]$-Poincar\'e lemma on Stein domains reduces these cohomologies to the holomorphic Koszul cohomology of $\theta$.
Third, a sheaf-theoretic computation localizes the Koszul cohomology at the critical points of $\theta$ and shows that it is concentrated in degree $n$ and finite-dimensional.
The strong Hodge decomposition then follows from the weak one \eqref{weak_Hodge}.

We will use the following tameness condition which was introduced by Dai and Yan in \parencite{DaiYan2022}. 

\begin{definition}\label{def:well-tamed}
    A Novikov type LG model $(M,g,\theta)$ is called well-tamed if there exists a compact subset $K \subset M$ and positive constants $C_1, C_2$ such that for every $x \in M \setminus K$ the following two inequalities hold:
    \begin{enumerate}
        \item $|\theta(x)| > C_1$, and
        \item $|\theta(x)|^2 - C_2\,|\nabla \theta(x)| > 0$,
    \end{enumerate}
    where $|\theta|$ and $|\nabla \theta|$ are taken with respect to the K\"ahler metric $g$.
\end{definition}

We adopt the definition of Stein manifolds from \parencite[Chapter~IV, Definition~1]{GrauertRemmert2004}, i.e complex manifolds that satisfy Cartan's Theorem~B.
Now we can state the main theorem.

\begin{theorem}\label{main_theorem_1}
	Let $(M,g,\theta)$ be a well-tamed Novikov type LG model. 
	Assume that $M$ is a Stein manifold. 
	Then we have the strong Hodge decomposition
	\begin{equation}
		L^2(\Lambda^\bullet(M)) = \H_\theta^\bullet \oplus \mathrm{im}(\tcr[\theta]) \oplus \mathrm{im}(\tcrd[\theta]).
	\end{equation}
	In particular, we have the isomorphism
	\begin{equation}
		\H_\theta^\bullet \cong \begin{cases}
			0, & k \ne n, \\
			\bigoplus_{x \in \mathrm{Crit}(\theta)} \mathcal{R}_{\theta,x}, & k=n.
		\end{cases}.
	\end{equation}
\end{theorem}

\subsection{Cohomological reduction}

We denote by $\Omega^\bullet(M)$ the space of smooth forms, by $\Omega^\bullet_0(M)$ the space of smooth forms with compact support, and by $\Omega^\bullet_{(2)}(M)=\Omega^\bullet(M) \cap L^2(\Lambda^\bullet(M))$ the space of smooth $L^2$-forms. 
We consider the following cohomologies associated with $\tcr[\theta]$.
\begin{gather}
	H^\bullet_{\tcr[\theta]} (M) = H^\bullet (\Omega^\bullet (M), \tcr[\theta]), \\
	H^\bullet_{(0, \tcr[\theta])} (M) = H^\bullet (\Omega^\bullet_0 (M), \tcr[\theta]), \\
	H^\bullet_{((2), \tcr[\theta])} (M) = H^\bullet (\Omega^\bullet_{(2)} (M), \tcr[\theta]).
\end{gather}
Moreover, we consider the following $L^2$-complex
\begin{equation}
	\cdots \to L^2 (\Lambda^{k-1} (M)) \overset{\tcr[\theta]}{\to} L^2 (\Lambda^{k} (M)) \overset{\tcr[\theta]}{\to} L^2 (\Lambda^{k+1} (M)) \to \cdots.
\end{equation}
We define its cohomology by $H^\bullet_{((2),\tcr[\theta],\#)} (M) = \ker (\tcr[\theta]) / \mathrm{im} (\tcr[\theta])$.

\begin{proposition}
	We have the isomorphism of cohomologies:
	\begin{equation}
		H^\bullet_{((2),\tcr[\theta])} (M) \cong H^\bullet_{((2),\tcr[\theta],\#)} (M).
	\end{equation}
\end{proposition}

\begin{proof}
	See \parencite[Lemma~2.51]{fan2011}.
\end{proof}

Next, we establish isomorphisms between the three types of $\tcr[\theta]$-cohomologies.
We use a homotopy argument, which generalizes the homotopy argument in \parencite{FanFang2021} for the original LG model. 
We first establish the following homotopy lemma.

\begin{lemma}\label{homotopy}
	Let $\rho$ be a smooth function with compact support in $M$ which equals $1$ in a neighborhood of $\mathrm{Crit}(\theta)$. 
	We define the operators
	\begin{gather}
		V_\theta = \frac{\iota_{\bar\theta^\sharp}}{|\theta|^2} = \frac{1}{|\theta|^2} \bar\theta_j g^{\bar jk} \ctr{k}, \\
		T_{\rho, \theta} = \rho + (\bar\partial \rho) V_\theta \frac{1}{1 + [\bar\partial, V_\theta]}, \\
		R_{\rho, \theta} = (1 - \rho) V_\theta \frac{1}{1 + [\bar\partial, V_\theta]}.
	\end{gather}
	Then we have
	\begin{gather}
		[\theta\wedge, V_\theta] = 1, \\
		[\bar\partial, [\bar\partial, V_\theta]] = [\theta\wedge, [\bar\partial, V_\theta]] =[V_\theta, [\bar\partial, V_\theta]] = 0, \\
		[\tcr[\theta], R_{\rho, \theta}] = 1 - T_{\rho, \theta}.
	\end{gather}
\end{lemma}

\begin{proof}
	We compute directly
	\begin{gather}
		[\theta \wedge, V_\theta] = [\theta_i \d z^i \wedge, \frac{1}{|\theta|^2} \bar\theta_j g^{\bar jk} \ctr{k}] = \frac{1}{|\theta|^2} \theta_k \bar\theta_j g^{\bar jk} = 1, \\
		[\bar\partial, V_\theta] = [\nabla_{\bar i} \d z^{\bar i} \wedge, \frac{1}{|\theta|^2} \bar\theta_j g^{\bar jk} \ctr{k}] = g^{\bar jk} (\nabla_{\bar i} \frac{\bar\theta_j}{|\theta|^2}) \d z^{\bar i} \wedge \ctr{k}.
	\end{gather}
	We then obtain
	\begin{gather}
		[\bar\partial, [\bar\partial, V_\theta]] = [\nabla_{\bar \ell} \d z^{\bar \ell} \wedge, g^{\bar jk} (\nabla_{\bar i} \frac{\bar\theta_j}{|\theta|^2}) \d z^{\bar i} \wedge \ctr{k}] = g^{\bar jk} (\nabla_{\bar \ell} \nabla_{\bar i} \frac{\bar\theta_j}{|\theta|^2} ) \d z^{\bar \ell} \wedge \d z^{\bar i} \wedge \ctr{k} = 0, \\
		[\theta \wedge, [\bar\partial, V_\theta]] = [\theta_\ell  \d z^\ell, g^{\bar jk} (\nabla_{\bar i} \frac{\bar\theta_j}{|\theta|^2}) \d z^{\bar i} \wedge \ctr{k}] = - \theta_k g^{\bar jk} (\nabla_{\bar i} \frac{\bar\theta_j}{|\theta|^2}) \d z^{\bar i} \wedge = 0, \\
		[V_\theta, [\bar\partial, V_\theta]] = [\frac{1}{|\theta|^2} \bar\theta_\ell  g^{\bar \ell p} \ctr{p}, g^{\bar jk} (\nabla_{\bar i} \frac{\bar\theta_j}{|\theta|^2}) \d z^{\bar i} \wedge \ctr{k}] = 0.
	\end{gather}
	We further compute
	\begin{equation}
		\begin{split}
			[\tcr[\theta], R_{\rho, \theta}] 
			& = [\tcr[\theta], 1-\rho] V_\theta \frac{1}{1 + [\bar\partial, V_\theta]} + (1 - \rho) [\tcr[\theta], V_\theta] \frac{1}{1 + [\bar\partial, V_\theta]} \\
			& = (- \bar\partial \rho) V_\theta \frac{1}{1 + [\bar\partial, V_\theta]} + (1 - \rho) (1 + [\bar\partial, V_\theta]) \frac{1}{1 + [\bar\partial, V_\theta]} \\
			& = 1 - T_{\rho, \theta}.
		\end{split}
	\end{equation}
\end{proof}

\begin{proposition}\label{isomorphism_cohomologies}
	Let $(M,g,\theta)$ be a well-tamed Novikov type LG model. 
	Then we have the isomorphisms
	\begin{equation}
		H^\bullet_{((2),\tcr[\theta])} (M) \cong H^\bullet_{\tcr[\theta]} (M) \cong H^\bullet_{(0,\tcr[\theta])} (M).
	\end{equation}
\end{proposition}

\begin{proof}
	The Lemma~\ref{homotopy} shows that $T_{\rho, \theta}: (\Omega^\bullet (M), \tcr[\theta]) \to (\Omega^\bullet_{0} (M), \tcr[\theta])$ is a chain map, and the operator $R_{\rho, \theta}$ defines a homotopy from $\mathrm{Id}_{\Omega^\bullet (M)}$ to $j \circ T_{\rho, \theta}$, where $j$ is the inclusion. 
	It follows that
	\begin{equation}
		H^\bullet_{\tcr[\theta]} (M) \cong H^\bullet_{(0,\tcr[\theta])} (M).
	\end{equation}
	It remains to show that $R_{\rho, \theta}$ is bounded. 
	Since $[\bar\partial, V_\theta]^n = 0$, we have
	\begin{equation}
		R_{\rho, \theta} = (1 - \rho) V_\theta \sum_{i=0}^{n-1} (-1)^i [\bar\partial, V_\theta]^i.
	\end{equation}
	It suffices to prove that $V_\theta$ and $[\bar\partial, V_\theta]$ are bounded. 
	Since
	\begin{equation}
		|[\bar\partial, V_\theta]| \le 3 \frac{|\nabla\theta|}{|\theta|^2},
	\end{equation}
	and the model is well-tamed, there exists $C>0$ such that on $\mathrm{supp}(1-\rho)$ we have
	\begin{equation}
		\frac{1}{|\theta|} < C,\ \frac{|\nabla\theta|}{|\theta|^2} < C,
	\end{equation}
	we obtain
	\begin{equation}
		|V_\theta| < C,\quad |[\bar\partial, V_\theta]| < 3C.
	\end{equation}
	Hence $R_{\rho, \theta}$ is bounded. 
	This completes the proof.
\end{proof}

Finally, we can reduce the $\tcr[\theta]$-cohomology to the Koszul cohomology. 

\begin{proposition}[{$\tcr[\theta]$} Poincar\'e lemma]\label{Omega_Poincare}
	Assume that $U \subset M$ is a Stein domain. 
    Let $\phi$ be a $\tcr[\theta]$-closed $k$-form on $U$. 
    Then there exist a $(k-1)$-form $\psi$ and a holomorphic $k$-form $\alpha$ such that
	\begin{equation}
		\phi = \alpha + \tcr[\theta] \psi.
	\end{equation}
\end{proposition}

\begin{proof}
	We write $\phi$ in the form
	\begin{equation}
		\phi = \sum_{p=0}^k \phi^{p,k-p},\ \phi^{p,k-p} \in \Omega^{p,k-p} (U).
	\end{equation}
	Then $\tcr[\theta] \phi = 0$ is equivalent to
	\begin{equation}
		\begin{cases}
			\bar\partial \phi^{0,k} = 0, \\
			\bar\partial \phi^{p,k-p} + \theta \wedge  \phi^{p-1,k-p+1} = 0,\ 1\le p \le k-1.
		\end{cases}
	\end{equation}
	Since $U$ is Stein, there exists a $(0,k-1)$-form $\psi^{0,k-1}$ such that $\phi^{0,k} = \bar\partial \psi^{0,k-1}$.
	Then we obtain
	\begin{equation}
		0 = \bar\partial \phi^{1,k-1} + \theta \wedge \bar\partial \psi^{0,k-1} = \bar\partial( \phi^{1,k-1} - \theta \wedge \psi^{0,k-1} ).
	\end{equation}
	Since $U$ is Stein, there exists a $(1,k-2)$-form $\psi^{1,k-2}$ such that
	\begin{equation}
		\phi^{1,k-1} = \bar\partial \psi^{1,k-2} + \theta \wedge \psi^{0,k-1}.
	\end{equation}
	Repeating this procedure, there exists a family of forms $\{\psi^{p,k-p-1}\}_{p=1}^{k-1}$, where $\psi^{p,k-p-1}$ is a $(p,k-p-1)$-form, such that
	\begin{equation}
		\phi^{p,k-p} = \bar\partial \psi^{p,k-p-1} + \theta \wedge \psi^{p-1,k-p}.
	\end{equation}
	Finally, we obtain
	\begin{equation}
		\bar\partial ( \phi^{k,0} - \theta \wedge \psi^{k-1,0} ) = 0.
	\end{equation}
	We set $\alpha = \phi^{k,0} - \theta \wedge \psi^{k-1,0}$, which is holomorphic, and $\psi = \sum_{p=0}^{k} \psi^{p,k-p-1}$.
	Then we have
	\begin{equation}
		\phi = \alpha + \tcr[\theta] \psi.
	\end{equation}
\end{proof}

\begin{proposition}\label{Koszul_reduction}
    Let $(M,g,\theta)$ be a Novikov type LG model, and $M$ be a Stein manifold.
    Then we have the isomorphism
    \begin{equation}
		H^\bullet_{\tcr[\theta]}(M)
		\cong
		H^\bullet(K^\bullet(\theta,M)),
	\end{equation}
	where
	\begin{equation}
		K^\bullet(\theta,M)
		=
		(\Omega^\bullet(M),\,\theta\wedge)
	\end{equation}
	is the holomorphic Koszul complex.
\end{proposition}

\begin{proof}
    Since $M$ is Stein, the $\tcr[\theta]$-Poincar\'e lemma (Proposition~\ref{Omega_Poincare}) implies that every $\tcr[\theta]$-cohomology class admits a representative of the form
	\begin{equation}
		\phi=\alpha+\tcr[\theta]\psi,
	\end{equation}
	where $\alpha$ is holomorphic. 
    Since $\tcr[\theta]\phi=0$, we obtain
	\begin{equation}
		0=\tcr[\theta]\alpha
		=\bar\partial\alpha+\theta\wedge\alpha.
	\end{equation}
	Since $\alpha$ is holomorphic, $\bar\partial\alpha=0$, hence $\theta\wedge\alpha=0$. 
	It follows that replacing $\alpha$ by another holomorphic representative changes it by an element of $\theta\wedge\Omega^{k-1}(M)$. 
    Then the correspondence $[\phi] \mapsto [\alpha]$ defines the isomorphism between $H^\bullet_{\tcr[\theta]}(M)$ and $H^\bullet(K^\bullet(\theta,M))$.
\end{proof}

\subsection{Finite-dimensionality and strong Hodge decomposition}

We can compute the Koszul cohomology of $\theta$ by localizing at critical points.

\begin{proposition}\label{finite_critical}
	Let $(M,g,\theta)$ be a well-tamed Novikov type LG model. 
	Assume that $M$ is a Stein manifold. 
	Then $\theta$ has only finitely many critical points.
\end{proposition}

\begin{proof}
	We denote by $S \subset M$ the analytic subset consisting of the critical points of $\theta$. 
	Since $\theta$ is well-tamed, $S$ is contained in a compact subset of $M$. 
	Since $S$ is closed, it is therefore compact. 
	By \parencite[Chapter~V, Theorem~3]{GrauertRemmert2004}, every compact analytic subset of a Stein manifold is finite. 
	Hence $S$ is finite.
\end{proof}

\begin{proposition}\label{Koszul_cohomology}
	Assume that $M$ is a Stein manifold of complex dimension $n$, and that the holomorphic $1$-form $\theta$ has finitely many critical points. 
    Then there is a noncanonical isomorphism of vector spaces
	\begin{equation}\label{Koszul_main_iso}
		H^k(K^\bullet(\theta,M))\cong
		\begin{cases}
			0,                                        & k\ne n, \\
			\displaystyle\bigoplus_{x\in\mathrm{Crit}(\theta)}
			\mathcal{R}_{\theta,x}, & k=n.
		\end{cases}
	\end{equation}
\end{proposition}

\begin{proof}
	We denote by $\mathcal{K}_\theta^\bullet$ the corresponding complex of sheaves,
	\begin{equation}
		\mathcal{K}_\theta^\bullet=(\Omega^\bullet,\theta\wedge).
	\end{equation}

	Since $M$ is Stein, the cohomology groups of coherent sheaves on $M$ vanish in positive degree. 
	Let $\{U_i\}$ be a Stein open covering of $M$. 
    We consider the \v{C}ech double complex associated to $\mathcal{K}_\theta^\bullet$ with respect to this covering. 
	Since each $\mathcal{K}^k$ and $\mathcal{H}^k(\mathcal{K}_\theta^\bullet)$ are coherent, their higher \v{C}ech cohomology groups vanish. 
	Then the total cohomology of the double complex reduces to the Koszul cohomology of global sections or the global sections of the Koszul cohomology of sheaves. 
	Hence for any Stein open set $U$, we have
	\begin{equation}\label{cohomology_section_commutativity}
		H^\bullet(K^\bullet(\theta,U))
		\cong
		\Gamma\bigl(\mathcal{H}^\bullet(\mathcal{K}_\theta^\bullet),U\bigr).
	\end{equation}

	To compute the right-hand side, let $(U,(z^1,\dots,z^n))$ be a coordinate chart.
	Writing
	\begin{equation}
		\theta=\sum_{i=1}^n \theta_i\,dz^i,
	\end{equation}
	we denote by $\mathcal{J}_\theta(U)$ the ideal generated by $\theta_1,\ldots,\theta_n$ on $U$. 
	This defines a coherent ideal sheaf $\mathcal{J}_\theta$ on $M$. 

	Let $x\in M$ and let $U\subset M$ be a sufficiently small Stein coordinate neighbourhood of $x$, containing no critical point of $\theta$ other than $x$. 
	The de Rham division lemma \parencite[cf.][Section~4.3]{Kulikov_1998} implies
	\begin{equation}
		H^k(K^\bullet(\theta,U))\cong
		\begin{cases}
			0,                                      & k<n, \\
			\mathcal{O}_M(U)/\mathcal{J}_\theta(U), & k=n.
		\end{cases}
	\end{equation}
	By \eqref{cohomology_section_commutativity}, the cohomology sheaves vanish for $k<n$ and are locally isomorphic to $\mathcal{R}_\theta$ for $k=n$. 
	Since these sheaves are supported on finitely many points, the local identifications give the following noncanonical isomorphisms
	\begin{equation}
		\mathcal{H}^k(\mathcal{K}_\theta^\bullet)\cong
		\begin{cases}
			0,                                & k<n, \\
			\mathcal{R}_\theta, & k=n.
		\end{cases}
	\end{equation}

	Taking global sections, we obtain
	\begin{equation}
		\Gamma(M,\mathcal{R}_\theta)
		\cong
		\bigoplus_{x\in\mathrm{Crit}(\theta)}
		\mathcal{R}_{\theta,x}.
	\end{equation}
	Combining these identifications yields \eqref{Koszul_main_iso}.
\end{proof}

Now we can complete the proof of Theorem~\ref{main_theorem_1}.

\begin{proof}[Proof of the main theorem]
	Proposition \ref{finite_critical} implies that $\theta$ has only finitely many critical points. 
	By Proposition \ref{Koszul_cohomology}, $H^\bullet(K^\bullet(\theta,M))$ is finite-dimensional. 
    By Proposition \ref{Koszul_reduction},    
    $H^\bullet_{\tcr[\theta]} (M)$ is finite-dimensional. 
	By Proposition \ref{isomorphism_cohomologies}, $H^\bullet_{((2), \tcr[\theta])} (M)$ is also finite-dimensional. 
	It follows that $\mathrm{im}(\tcr[\theta])$ is closed, and similarly $\mathrm{im}(\tcrd[\theta])$ is closed. 
	Therefore, the strong Hodge decomposition follows from \eqref{weak_Hodge}. 
	This yields $\ker(\tcr[\theta]) = \H^\bullet_\theta \oplus \mathrm{im}(\tcr[\theta])$, and hence $\H^\bullet_\theta \cong H^\bullet_{((2),\tcr[\theta])} (M)$.
\end{proof}

\section{Hyperplane Arrangements, Complete K\"ahler Metrics and Well-Tameness}\label{sec:YY}

We first recall the intersection poset and the critical-point count for logarithmic master functions.
We then construct a complete K\"ahler metric on an arrangement complement.
For a simple normal crossing arrangement, weighted estimates in a logarithmic frame give a uniform bound for the covariant derivative of the master one-form and a positive lower bound for its norm outside a compact set.
These two estimates establish well-tameness in the sense of Definition~\ref{def:well-tamed}.
The metric has cylindrical contributions in normal directions, and its tangential coefficients grow near the boundary; the proof retains this growth throughout.

\subsection{The intersection poset and the M\"obius function}\

Let homogeneous coordinates be
\begin{equation*}
	[Z]=[Z_0:\cdots:Z_n]\in \P^n,
\end{equation*}
and let each hyperplane be
\begin{equation}
	H_i=\{f_i(Z)=0\},\qquad
	f_i(Z)=\sum_{j=0}^n c_{ij}Z_j .
\end{equation}

Then the set $\mathcal{A}=\{H_1,\dots,H_m\}$ is a hyperplane arrangement in $\P^n$.
Its intersection poset $L(\mathcal{A})$ consists of all nonempty intersections of hyperplanes:
\begin{equation}
	X=\bigcap_{i\in I} H_i,
\end{equation}
together with the ambient space $\P^n$, ordered by reverse inclusion:
\begin{equation*}
	X\leq Y \quad \Longleftrightarrow \quad Y\subseteq X.
\end{equation*}
So the minimal element is $\hat{0}=\P^n$.
The M\"obius function $	\mu:L(\mathcal{A})\to \Z$ is defined recursively by $\mu(\P^n)=1$ and, for $X>\P^n$,
\begin{equation*}
	\mu(X)
	=
	-\sum_{Y<X}\mu(Y).
\end{equation*}
Equivalently,
\begin{equation*}
	\sum_{Y\leq X}\mu(Y)=0
	\qquad\text{for }X\neq \P^n.
\end{equation*}
More generally, the two-variable M\"obius function of the poset is $\mu(Y,X)$, defined for $Y\leq X$ by $\mu(X,X)=1$, and
\begin{equation*}
	\mu(Y,X)
	=
	-\sum_{Y\leq Z<X}\mu(Y,Z).
\end{equation*}
The one-variable version above is
\begin{equation*}
	\mu(X)=\mu(\P^n,X).
\end{equation*}

For example:
\begin{enumerate}
	\item $\mu(\P^n)=1$.
	\item If $H_i$ is a hyperplane, then
	$\mu(H_i)=-1$.
	\item If $X$ has codimension $2$ and is contained in exactly $r$ hyperplanes, then
	\begin{equation*}
		\mu(X)=-(1-r)=r-1.
	\end{equation*}
\end{enumerate}
For a systematic account, see \parencite{OrlikTerao1992}.

\subsection{Master functions and critical-point counts}

For weights $a_r\in\C^*$, consider the multivalued master function
\begin{equation}
 F=\sum_{r=1}^m a_r\log f_r,
 \qquad \theta=\d F=\sum_{r=1}^m a_r\frac{\d f_r}{f_r}
\end{equation}
on the complement of the corresponding central arrangement in $\C^{n+1}$.
The one-form $\theta$ is invariant under scalar multiplication.
Its contraction with the Euler vector field is $\sum_r a_r$, so it descends to
$M=\P^n\setminus\bigcup_rH_r$ if and only if
\begin{equation}\label{eq:sum-weights}
 \sum_{r=1}^m a_r=0.
\end{equation}
We impose this condition throughout.
Choose $H_m$ as the hyperplane at infinity.
Then $M$ is a principal affine open subset of $\C^n$, hence is Stein.

\begin{definition}[Nonresonant weights]\label{subsec:resonance}
The weights are called nonresonant in this section if, for every nonempty proper intersection $X\ne\P^n$ of the arrangement,
\begin{equation}\label{eq:nonresonance-weights}
 a_X:=\sum_{H_r\supset X}a_r\ne0.
\end{equation}
Otherwise they are called resonant.
\end{definition}

We apply the Varchenko--Orlik--Terao critical-point counting formula.
With $H_m$ chosen as the hyperplane at infinity, its weight is
\begin{equation*}
 a_m=-\sum_{r=1}^{m-1}a_r.
\end{equation*}
The nonresonance condition \eqref{eq:nonresonance-weights} implies that every dense edge of the projective arrangement has nonzero weight, so the resulting affine weighted arrangement is \emph{unbalanced} in Varchenko's terminology.
It is also essential: a common nonempty intersection of all projective hyperplanes would have weight $\sum_ra_r=0$, contradicting nonresonance.
Thus, by \parencite[Lemma~2.8]{Varchenko2015CriticalSet},
the critical set $\mathrm{Crit}(\theta)$ is finite and
\begin{equation}\label{eq:logarithmic-count}
 \mu_\theta=\sum_{p\in \mathrm{Crit}(\theta)}\mu_p=(-1)^n\chi(M),
\end{equation}
where $\mu_p$ is the local Milnor number of $F$ at $p$.
In the elementary case of $n$ affine hyperplanes in essential position, affine coordinates identify $M$ with $(\C^*)^n$ and $\theta$ with $\sum_{r=1}^n a_r\,\d z_r/z_r$; its critical set is empty and both sides of \eqref{eq:logarithmic-count} vanish.
For sufficiently general weights, the critical points are nondegenerate by \parencite[Theorem~2.7]{Varchenko2015CriticalSet}, so the total Milnor number equals the number of distinct critical points.

For a simple normal crossing arrangement, the counting formula \eqref{eq:logarithmic-count} also holds whenever all $a_r\ne0$: the logarithmic differential has nonzero boundary residues on $(\P^n,D)$, and the logarithmic Poincar\'e--Hopf formula applies.
For this geometric form of the counting formula, see \parencite[Section~2]{Huh2013MaximumLikelihood}.

The Euler characteristic is computed by the intersection poset:
\begin{equation}
 \chi(M)=\sum_{X\in L(\mathcal A)}\mu(X)\chi(X),
 \qquad \chi(X)=\dim X+1.
\end{equation}
Here the sum includes $X=\P^n$.
For nonresonant weights, or for a simple normal crossing arrangement with nonzero weights,
\begin{equation}\label{eq:counting-lattice}
 \sum_{p\in \mathrm{Crit}(\theta)}\mu_p
 =(-1)^n\sum_{X\in L(\mathcal A)}\mu(X)(\dim X+1).
\end{equation}
The circle-valued Morse-theoretic approach of \parencite{KohnoPajitnov2015} is discussed in the Introduction.

\subsection{A complete K\"ahler metric}\label{sec:cylindrical}\

We construct a complete K\"ahler metric $g$ on the complement of an arbitrary hyperplane arrangement using a potential quadratic in logarithms of section norms.

Let $H_r=\{f_r=0\}\subset\P^n$, $r=1,\dots,m$, where
\begin{equation*}
 f_r(Z)=\sum_{j=0}^n c_{rj}Z_j.
\end{equation*}
Give $\O(1)$ the Fubini--Study Hermitian metric.
Denote its curvature form by $\omega_{FS}$ and the associated metric by $g_{FS}$.
Write
\begin{equation}
 \|f_r([Z])\|^2_{FS}=\frac{|f_r(Z)|^2}{\|Z\|^2},
\end{equation}

Write
\begin{equation}
 D=\bigcup_{r=1}^mH_r,\qquad M=\P^n\setminus D
\end{equation}
Set
\begin{equation*}
 C_r=1+\log\left(\sum_{j=0}^n|c_{rj}|^2\right),
\end{equation*}
and define
\begin{equation}\label{eq:logarithmic-rho}
 \rho_r=C_r-\log\|f_r\|^2_{FS}
 =1-\log\frac{|f_r(Z)|^2}
                  {\left(\sum_{j=0}^n|c_{rj}|^2\right)\|Z\|^2}.
\end{equation}
\begin{equation*}
 |f_r(Z)|^2\le
 \left(\sum_{j=0}^n|c_{rj}|^2\right)\|Z\|^2,
\end{equation*}
so $\rho_r\ge1$ on $M$.
This definition is unchanged when $f_r$ is multiplied by a nonzero constant.
Set
\begin{equation}\label{eq:cylindrical-form}
 \omega=\omega_{FS}+\frac{\i}{2}\sum_{r=1}^m
                 \partial\bar\partial(\rho_r^2).
\end{equation}
The product rule gives
\begin{equation}
 \frac{\i}{2}\partial\bar\partial(\rho_r^2)
 =\i\partial\rho_r\wedge\bar\partial\rho_r
  +\i\rho_r\partial\bar\partial\rho_r.
\end{equation}
On $M$, the curvature identity gives $\i\partial\bar\partial\rho_r=\omega_{FS}$.

Set $Q=\sum_r\rho_r$, so $Q\ge m$.
Then
\begin{equation}\label{eq:cylindrical-expanded}
 \omega=(1+Q)\omega_{FS}
          +\i\sum_{r=1}^m\partial\rho_r\wedge\bar\partial\rho_r.
\end{equation}
Thus $\omega$ is positive and closed, hence K\"ahler.
Let $g$ be its associated metric, regarded as a Hermitian metric on $T_\C M=TM\otimes_\R\C$.
For $v,w\in T^{1,0}M$,
\begin{equation}\label{eq:metric-Q}
 g(v,w)=(1+Q)\,g_{FS}(v,w)
        +\sum_r\partial\rho_r(v)\,
                    \overline{\partial\rho_r(w)}.
\end{equation}
This determines $g$ on $T_\C M$ by orthogonality and conjugation.

To prove completeness, note that $\rho_r\to+\infty$ as a point approaches $H_r$, while all the remaining functions $\rho_r$ are bounded below by $1$.
Thus $Q\to+\infty$ as a point approaches the boundary $D$.
Since $\P^n$ is compact, every curve escaping all compact subsets of $M$ approaches $D$.

For a piecewise $C^1$ curve $\gamma:[0,T)\to M$ escaping every compact subset, $\frac{\d}{\d t}\bigl(Q(\gamma(t))\bigr)=\sum_r\dot\gamma\rho_r$.
Writing $\dot\gamma=v+\bar v$ with $v\in T^{1,0}M$, we have $\dot\gamma\rho_r=2\operatorname{Re}\partial\rho_r(v)$ and $g(\dot\gamma,\dot\gamma)=2g(v,v)$.
Together with \eqref{eq:metric-Q}, these identities give
\begin{equation*}
 |\frac{\d}{\d t}\bigl(Q(\gamma(t))\bigr)|^2
 \le m\sum_r|\dot\gamma\rho_r|^2
 \le2m\,g(\dot\gamma,\dot\gamma).
\end{equation*}

Consequently, for $0<t<T$, the arc length satisfies
\begin{equation*}
 \begin{aligned}
 \int_0^t\sqrt{g(\dot\gamma(s),\dot\gamma(s))}\,\d s
 &\ge\frac1{\sqrt{2m}}\int_0^t
       \left|\frac{\d}{\d s}\bigl(Q(\gamma(s))\bigr)\right|\,\d s\\
 &\ge\frac1{\sqrt{2m}}
       |Q(\gamma(t))-Q(\gamma(0))|.
 \end{aligned}
\end{equation*}
The boundary behavior gives $Q(\gamma(t))\to+\infty$ as $t\to T$, so every such curve has infinite arc length.
Hence $g$ is complete.

\subsubsection*{Explicit formulas in homogeneous coordinates}\label{sec5.2-expl-metr}

Let $\pi:\C^{n+1}\setminus\{0\}\to\P^n$ be the standard projection, and use the Hermitian inner product $\langle V,W\rangle=\sum_{j=0}^nV_j\overline{W_j}$ on $\C^{n+1}$.
For $Z$ with $[Z]\in M$, set $S=\langle Z,Z\rangle$; the definition of $\rho_r$ gives
\begin{equation*}
 \rho_r([Z])=\log S-\log|f_r(Z)|^2+C_r,\qquad
 (\pi^*\partial\rho_r)_Z(V)
 =\frac{\langle V,Z\rangle}{S}-\frac{f_r(V)}{f_r(Z)}.
\end{equation*}
For $v=d\pi_Z(V)$ and $w=d\pi_Z(W)$, equation~\eqref{eq:metric-Q} therefore becomes
\begin{equation}
 \begin{aligned}
 g(v,w)={}&(1+Q([Z]))
 \left(\frac{\langle V,W\rangle}{S}
       -\frac{\langle V,Z\rangle\langle Z,W\rangle}{S^2}\right)\\
 &+\sum_r
 \left(\frac{\langle V,Z\rangle}{S}-\frac{f_r(V)}{f_r(Z)}\right)
 \overline{\left(\frac{\langle W,Z\rangle}{S}
                         -\frac{f_r(W)}{f_r(Z)}\right)}.
 \end{aligned}
\end{equation}
The expression is unchanged when $V$ or $W$ is replaced by another lift differing by a multiple of $Z$.

On the affine chart $Z_0\ne0$, write
\begin{equation*}
 x_i=Z_i/Z_0,\qquad i=1,\dots,n.
\end{equation*}
We use the same symbol for the local representative
\begin{equation*}
 f_r(x):=f_r(1,x)=c_{r0}+\sum_{i=1}^n c_{ri}x_i,\qquad
 S=1+\sum_{i=1}^n|x_i|^2.
\end{equation*}
Then
\begin{equation*}
 \rho_r=\log S-\log|f_r(x)|^2+C_r,\qquad
 \partial_i\rho_r=\frac{\overline{x_i}}{S}-\frac{c_{ri}}{f_r(x)}.
\end{equation*}
With $g_{i\bar j}=g(\partial_{x_i},\partial_{x_j})$, the metric components are
\begin{equation}
 g_{i\bar j}
 =(1+Q)\left(\frac{\delta_{ij}}{S}
                -\frac{\overline{x_i}x_j}{S^2}\right)
  +\sum_r\left(\frac{\overline{x_i}}{S}-\frac{c_{ri}}{f_r}\right)
          \left(\frac{x_j}{S}
                   -\frac{\overline{c_{rj}}}{\overline{f_r}}\right).
\end{equation}
The logarithmic terms $\partial f_r/f_r$ give cylindrical contributions in normal directions, while the factor $(1+Q)$ multiplying $g_{FS}$ grows near the boundary.
The estimates below retain both contributions.

\subsection{Well-tameness in the normal-crossing case}\label{sec5.3}\label{subsec:norm-estimates}

Throughout this subsection, $D$ is a simple normal crossing hyperplane arrangement.
We use the metric of \eqref{eq:cylindrical-form} and write $\theta=\d F$.
All estimate constants below are independent of the point and of the tangent vectors; they may depend on the fixed arrangement, weights, metric and coordinate neighborhood.

Fix $p_0\in D$, and suppose that precisely $H_1,\ldots,H_k$ contain $p_0$.
On a sufficiently small coordinate neighborhood, choose
\begin{equation*}
 D=\{z_1\cdots z_k=0\},\qquad
 e_i=z_i\partial_{z_i}\ (i\le k),\qquad
 e_j=\partial_{z_j}\ (j>k).
\end{equation*}
The dual holomorphic coframe on the complement is
$\vartheta^i=\d z_i/z_i$ and $\vartheta^j=\d z_j$.
Work on a smaller neighborhood with compact closure in this coordinate chart, and choose a nonvanishing holomorphic frame $\tau$ of $\O(1)$.
Since $H_i=\{z_i=0\}$, we can write $f_i=z_iu_i\tau$ with $u_i$ holomorphic and nowhere vanishing, so $\rho_i=-\log|z_i|^2+\psi_i$, where $\psi_i=C_i-\log|u_i|^2-\log\|\tau\|_{FS}^2$ is smooth and real-valued.
For $r>k$, set $\psi_r=\rho_r$, which is smooth since $f_r$ is nonvanishing on this neighborhood.

\begin{lemma}[Weighted comparison of the metric]\label{lem:weighted-metric}
Define the positive scales by
\begin{equation}\label{eq:weighted-scales}
 \sigma_i^2=1+(1+Q)|z_i|^2\quad(i\le k),\qquad
 \sigma_j^2=1+Q\quad(j>k).
\end{equation}
There exist constants $c_g,C_g>0$, independent of the point in the coordinate neighborhood intersected with $M$ and of $v$, such that, for every $v=\sum_i v_i e_i$,
\begin{equation}\label{eq:weighted-comparison}
 c_g\sum_i\sigma_i^2|v_i|^2
 \le g(v,v)\le C_g\sum_i\sigma_i^2|v_i|^2.
\end{equation}
\end{lemma}

\begin{proof}
Put
\begin{equation*}
 B(v)=\sum_{i\le k}|z_i|^2|v_i|^2+\sum_{j>k}|v_j|^2.
\end{equation*}
This is the squared Euclidean norm of $v$ in the ordinary coordinates $z$.
The smooth positive definite metric $g_{FS}$ satisfies
\begin{equation*}
 c_{FS}B(v)\le g_{FS}(v,v)\le C_{FS}B(v)
\end{equation*}
for constants $c_{FS},C_{FS}>0$ uniform on the compact closure.
For a uniform constant $C_{\partial\psi}>0$, the bounded ordinary coordinate derivatives of $\psi_r$ and Cauchy--Schwarz give $|\partial\psi_r(v)|\le C_{\partial\psi}B(v)^{1/2}$.
Thus $C_\psi:=mC_{\partial\psi}^2$ satisfies
\begin{equation*}
 \sum_{r=1}^m|\partial\psi_r(v)|^2\le C_\psi B(v).
\end{equation*}

The identities
\begin{equation*}
 \partial\rho_i(v)=-v_i+\partial\psi_i(v)\quad(i\le k),
 \qquad \partial\rho_r(v)=\partial\psi_r(v)\quad(r>k)
\end{equation*}
therefore imply
\begin{equation*}
 \begin{aligned}
 \sum_{r=1}^m|\partial\rho_r(v)|^2
 &\le2\sum_{i\le k}|v_i|^2+2C_\psi B(v),\\
 \sum_{i\le k}|v_i|^2
 &\le2\sum_{r=1}^m|\partial\rho_r(v)|^2+2C_\psi B(v).
 \end{aligned}
\end{equation*}
Since $1+Q\ge1$, these inequalities and \eqref{eq:metric-Q} give constants $c_g,C_g>0$, independent of the point and of $v$, with
\begin{equation*}
 c_g\left((1+Q)B(v)+\sum_{i\le k}|v_i|^2\right)
 \le g(v,v)\le
 C_g\left((1+Q)B(v)+\sum_{i\le k}|v_i|^2\right).
\end{equation*}
Indeed, the first inequality for the logarithmic derivatives bounds $g(v,v)$ above by $(C_{FS}+2C_\psi)(1+Q)B(v)+2\sum_{i\le k}|v_i|^2$.
For the lower bound, $g(v,v)\ge c_{FS}(1+Q)B(v)$ and $g(v,v)\ge\sum_r|\partial\rho_r(v)|^2$ yield
\begin{equation*}
 (1+Q)B(v)+\sum_{i\le k}|v_i|^2
 \le\left(2+\frac{1+2C_\psi}{c_{FS}}\right)g(v,v).
\end{equation*}
Finally, the definition of the scales gives
\begin{equation*}
 (1+Q)B(v)+\sum_{i\le k}|v_i|^2
 =\sum_i\sigma_i^2|v_i|^2,
\end{equation*}
which proves the claim.
\end{proof}

\begin{lemma}[Connection estimates]\label{lem:weighted-connection}
Let $\nabla$ be the complexified Levi--Civita connection, whose restriction to $T^{1,0}M$ is the Chern connection.
In the frame above, define the coefficients by
\begin{equation*}
 G_{i\bar j}=g(e_i,e_j),\qquad
 \nabla_{e_i}e_j=\Gamma^p_{ij}e_p.
\end{equation*}
Let $c_g,C_g>0$ be comparison constants supplied by Lemma~\ref{lem:weighted-metric}.
There is a uniform constant $C_\Gamma>0$ such that
\begin{equation}\label{eq:weighted-inverse}
 |G^{\bar\ell p}|\le\frac{c_g^{-1}}{\sigma_\ell\sigma_p},
 \qquad
 |\Gamma^p_{ij}|\le C_\Gamma\frac{\sigma_i\sigma_j}{\sigma_p}.
\end{equation}
\end{lemma}

\begin{proof}
Let $W=\operatorname{diag}(\sigma_1,\ldots,\sigma_n)$ and $A=W^{-1}GW^{-1}$.
Applying Lemma~\ref{lem:weighted-metric} to $v=\sum_i(\xi_i/\sigma_i)e_i$ gives, for every $\xi\in\C^n$,
\begin{equation*}
 c_g\sum_i|\xi_i|^2
 \le\sum_{i,j}A_{i\bar j}\xi_i\overline{\xi_j}
 \le C_g\sum_i|\xi_i|^2.
\end{equation*}
The lower bound for $A$ gives
\begin{equation*}
 |A^{-1}\xi|\le c_g^{-1}|\xi|\quad(\xi\in\C^n),
 \qquad |(A^{-1})_{\ell p}|\le c_g^{-1},
\end{equation*}
where $|\cdot|$ is the Euclidean norm on $\C^n$.
Using $G^{-1}=W^{-1}A^{-1}W^{-1}$,
\begin{equation*}
 |G^{\bar\ell p}|
 =\frac{|(A^{-1})_{\ell p}|}{\sigma_\ell\sigma_p}
 \le\frac{c_g^{-1}}{\sigma_\ell\sigma_p}.
\end{equation*}

Choose uniform constants $C_{\partial\rho},C_{\partial Q}>0$ such that
\begin{equation*}
 |e_j\rho_r|+|e_i e_j\rho_r|+|e_i\bar e_\ell\rho_r|\le C_{\partial\rho},
 \qquad |e_iQ|\le C_{\partial Q}.
\end{equation*}
The first logarithmic derivatives of $-\log|z_i|^2$ for $i\le k$ are constant, and its second derivatives in these commuting holomorphic and antiholomorphic logarithmic frames vanish.
All remaining terms extend smoothly to the coordinate neighborhood.

Put $q_i=z_i$ for $i\le k$ and $q_i=1$ for $i>k$, and write $H_{j\bar\ell}=g_{FS}(e_j,e_\ell)$.
Since $e_i\bar q_\ell=0$,
\begin{equation*}
 \begin{aligned}
 H_{j\bar\ell}
 &=q_j\bar q_\ell g_{FS}(\partial_{z_j},\partial_{z_\ell}),\\
 e_iH_{j\bar\ell}
 &=(e_iq_j)\bar q_\ell g_{FS}(\partial_{z_j},\partial_{z_\ell})\\
 &\quad+q_iq_j\bar q_\ell\partial_{z_i}
             [g_{FS}(\partial_{z_j},\partial_{z_\ell})].
 \end{aligned}
\end{equation*}
The functions $q_i$, the ordinary-coordinate coefficients of $g_{FS}$, and their first derivatives are bounded on the compact closure.
Moreover, $e_iq_j=\delta_{ij}q_j$ for $j\le k$ and vanishes for $j>k$, so a uniform constant $C_H>0$ satisfies
\begin{equation*}
 |H_{j\bar\ell}|+|e_iH_{j\bar\ell}|
 \le C_H|q_jq_\ell|
 \le C_H\frac{\sigma_j\sigma_\ell}{1+Q}.
\end{equation*}
The last inequality uses $|q_i|\le\sigma_i/\sqrt{1+Q}$ for every $i$.
Now
\begin{equation*}
 G_{j\bar\ell}=(1+Q)H_{j\bar\ell}
       +\sum_r(e_j\rho_r)(\bar e_\ell\rho_r),
\end{equation*}
and therefore
\begin{align*}
 e_iG_{j\bar\ell}
 &=(e_iQ)H_{j\bar\ell}+(1+Q)e_iH_{j\bar\ell}\\
 &\quad+\sum_r\bigl[(e_ie_j\rho_r)(\bar e_\ell\rho_r)
                 +(e_j\rho_r)(e_i\bar e_\ell\rho_r)\bigr].
\end{align*}
\begin{equation}\label{eq:weighted-metric-derivative}
 \begin{aligned}
 |e_iG_{j\bar\ell}|
 &\le C_{\partial Q}C_H\frac{\sigma_j\sigma_\ell}{1+Q}
       +C_H\sigma_j\sigma_\ell+2mC_{\partial\rho}^2\\
 &\le C_{\partial G}\sigma_i\sigma_j\sigma_\ell,
 \end{aligned}
\end{equation}
where $C_{\partial G}:=C_{\partial Q}C_H+C_H+2mC_{\partial\rho}^2$, since $\sigma_i\ge1$ and $1+Q\ge1$.
The standard K\"ahler connection formula in the holomorphic frame gives
\begin{equation*}
 \Gamma^p_{ij}=(e_iG_{j\bar\ell})G^{\bar\ell p}.
\end{equation*}
\begin{equation*}
 \begin{aligned}
 |\Gamma^p_{ij}|
 &\le\sum_\ell|e_iG_{j\bar\ell}|\,|G^{\bar\ell p}|\\
 &\le C_{\partial G}c_g^{-1}\sum_\ell
       \frac{\sigma_i\sigma_j\sigma_\ell}{\sigma_\ell\sigma_p}
 =C_\Gamma\frac{\sigma_i\sigma_j}{\sigma_p},
 \end{aligned}
\end{equation*}
with $C_\Gamma:=nC_{\partial G}c_g^{-1}$.
\end{proof}

\begin{lemma}[Bounded covariant derivative]\label{lem:bounded-master-derivative}
The form $\theta$ has uniformly bounded covariant derivative in a neighborhood of $D$:
\begin{equation}\label{eq:nabla-decay}
 |\nabla\theta|_g\le C_{\nabla\theta}.
\end{equation}
\end{lemma}

\begin{proof}
Let $c_g$ and $C_\Gamma$ be the constants of Lemmas~\ref{lem:weighted-metric} and \ref{lem:weighted-connection}, respectively.
In a boundary coordinate neighborhood,
\begin{equation*}
 \theta=\sum_{i\le k}a_i\frac{\d z_i}{z_i}+\beta,
\end{equation*}
where $\beta$ is holomorphic across the boundary.
Writing $\theta=\sum_\ell\theta_\ell\vartheta^\ell$, choose a uniform constant $C_{\theta}>0$ with
$|\theta_\ell|+|e_j\theta_\ell|\le C_{\theta}$.
For $(\nabla\theta)_{j\ell}:=(\nabla_{e_j}\theta)(e_\ell)$, we have
\begin{equation*}
 \begin{aligned}
 (\nabla\theta)_{j\ell}&=e_j\theta_\ell-\Gamma^i_{j\ell}\theta_i,\\
 |(\nabla\theta)_{j\ell}|
 &\le C_e
       \left(1+\sigma_j\sigma_\ell\sum_i\sigma_i^{-1}\right),
 \end{aligned}
\end{equation*}
where $C_e:=C_{\theta}\max\{1,C_\Gamma\}$.
Since $\nabla_{\bar e_j}\theta=0$, the dual metric comparison from \eqref{eq:weighted-comparison} and its tensor square give
\begin{align*}
 |\nabla\theta|_g^2
 &\le c_g^{-2}\sum_{j,\ell}
       \frac{|(\nabla\theta)_{j\ell}|^2}{\sigma_j^2\sigma_\ell^2}\\
 &\le2c_g^{-2}C_e^2
       \sum_{j,\ell}\left[\frac1{\sigma_j^2\sigma_\ell^2}
                  +\left(\sum_i\frac1{\sigma_i}\right)^2\right].
\end{align*}
As $\sigma_i\ge1$, the last sum is at most $n^2(1+n^2)$.
Taking a finite boundary cover, choose $C_{\nabla\theta}>0$ so that, on each chart,
\begin{equation*}
 C_{\nabla\theta}^2\ge
 2c_g^{-2}C_e^2n^2(1+n^2).
\end{equation*}
\end{proof}

\begin{lemma}[Positive lower bound]\label{lem:master-lower-bound}
If every $a_r\ne0$, there exist $K\Subset M$ and $c_{\theta}>0$ such that
\begin{equation}\label{eq:norm-lower-bound}
 |\theta|_g^2\ge c_{\theta}\qquad\text{on }M\setminus K.
\end{equation}
\end{lemma}

\begin{proof}
In each boundary coordinate neighborhood, all the smooth functions $\psi_r$ are bounded.
At each point choose $p\le k$ such that $\rho_p=\max_{i\le k}\rho_i$.
Choose uniform constants $C_Q,C_z>0$ such that
\begin{equation*}
 \begin{gathered}
 1+Q\le C_Q(1+\rho_p),\qquad
 |z_p|^2=e^{\psi_p-\rho_p}\le C_ze^{-\rho_p},\\
 (1+Q)|z_p|^2\le C_QC_z(1+\rho_p)e^{-\rho_p}\longrightarrow0
 \end{gathered}
\end{equation*}
as the point approaches the boundary.
For this choice of $p$, the background metric has a uniform bound $C_{FS}>0$ with
\begin{equation*}
 H_{p\bar p}=g_{FS}(e_p,e_p)\le C_{FS}|z_p|^2.
\end{equation*}
Formula~\eqref{eq:metric-Q} and the bound $C_{\partial\rho}$ from the proof of Lemma~\ref{lem:weighted-connection} give
\begin{equation*}
 \begin{aligned}
 G_{p\bar p}
 &=(1+Q)H_{p\bar p}+\sum_r|e_p\rho_r|^2\\
 &\le C_{FS}C_QC_z(1+\rho_p)e^{-\rho_p}+mC_{\partial\rho}^2
 \le C_p,
 \end{aligned}
\end{equation*}
where $C_p:=C_{FS}C_QC_z+mC_{\partial\rho}^2$, using $(1+\rho_p)e^{-\rho_p}\le1$.
The local expression for $\theta$ also gives $|\theta(e_p)-a_p|\le C_\beta|z_p|$ for a uniform bound $C_\beta>0$ on the ordinary coefficients of $\beta$.
With $a_*:=\min_r|a_r|>0$, sufficiently near the boundary,
\begin{equation}\label{eq:norm-smooth}
 |\theta(e_p)|\ge\frac{|a_p|}{2}\ge\frac{a_*}{2}.
\end{equation}
\begin{equation}\label{eq:norm-stratum}
 |\theta|_g^2\ge\frac{|\theta(e_p)|^2}{G_{p\bar p}}
 \ge\frac{a_*^2}{4C_p}>0.
\end{equation}
Choose a finite collection of smaller boundary coordinate neighborhoods and take $c_{\theta}$ to be the minimum of their lower bounds $a_*^2/(4C_p)$.
The required thresholds can also be made uniform.
The complement of a sufficiently small boundary neighborhood is compact in $M$, proving the claim.
\end{proof}

\begin{proposition}\label{prop:well-tamed-normal-crossing}
Let $D=\bigcup_{r=1}^mH_r\subset\P^n$ be a simple normal crossing hyperplane arrangement, and let $a_r\in\C^*$ satisfy $\sum_ra_r=0$.
For the complete K\"ahler metric of \eqref{eq:cylindrical-form}, the triple $(M,g,\theta)$, where $M=\P^n\setminus D$ and $\theta=\sum_ra_r\partial f_r/f_r$, is a well-tamed Novikov type LG model in the sense of Definition~\ref{def:well-tamed}.
\end{proposition}

\begin{proof}
The form $\theta$ is closed and holomorphic on $M$ and is nontrivial because its residues are nonzero.
Completeness was proved above.
The two preceding lemmas give, outside a compact subset,
$|\theta|^2\ge c_{\theta}$ and $|\nabla\theta|\le C_{\nabla\theta}$.
Set $C_1=\sqrt{c_{\theta}}/2$ and $C_2=c_{\theta}/(2(C_{\nabla\theta}+1))$.
Then $|\theta|>C_1$ and
\begin{equation}\label{eq:well-tamed-estimate}
 |\theta|^2-C_2|\nabla\theta|
 \ge c_{\theta}-\frac{c_{\theta}C_{\nabla\theta}}{2(C_{\nabla\theta}+1)}\ge\frac{c_{\theta}}{2}>0,
\end{equation}
as required.
\end{proof}

\begin{corollary}
Under the hypotheses of Proposition~\ref{prop:well-tamed-normal-crossing},
\begin{equation*}
 \dim_\C\mathcal H^n_\theta=\sum_{p\in \mathrm{Crit}(\theta)}\mu_p=(-1)^n\chi(M),
 \qquad \mathcal H^j_\theta=0\quad(j\ne n).
\end{equation*}
\end{corollary}

\begin{proof}
Since $M$ is Stein, Proposition~\ref{prop:well-tamed-normal-crossing} allows us to apply Theorem~\ref{main_theorem_1}.
The affine volume form on $\P^n\setminus H_m\simeq\C^n$, restricted to $M$, trivializes the canonical bundle and gives the identification with the local Milnor rings.
The dimension formula follows from \eqref{eq:logarithmic-count}.
\end{proof}
\subsection{Non-normal-crossing case}\label{sec5.4}\

For a general arrangement, we pass to a log resolution and use a Poincar\'e type metric.
The curvature forms of exceptional divisor line bundles need not have the positivity used in Section~\ref{sec:cylindrical}.
The construction below controls these curvature terms by a smooth K\"ahler background on the resolution.

Let $D=\bigcup_{r=1}^m H_r\subset\P^n$ and choose a projective log resolution
\begin{equation*}
 \pi:\widetilde X\longrightarrow\P^n,\qquad
 E=(\pi^{-1}D)_{\mathrm{red}}=\bigcup_{\nu=1}^N E_\nu,
\end{equation*}
where $E$ is a simple normal crossing divisor and $\pi$ is an isomorphism over $M$.
Its irreducible components $E_\nu$ consist of the strict transforms of the hyperplanes $H_r$ and the exceptional divisors of $\pi$.
For example, one may use a De Concini--Procesi wonderful model \parencite{DeConciniProcesi1995Wonderful}.
We identify $M$ with $\widetilde X\setminus E$ and retain the notation $\theta=\d F$ for its pullback.
Each pullback divisor $\pi^*H_r$ is supported on $E$, and $\theta$ extends logarithmically across $E$.
Its residue along $E_\nu$ is the constant
\begin{equation}\label{eq:resolution-residues}
 \lambda_\nu=\sum_{r=1}^m a_r\,
                  \operatorname{ord}_{E_\nu}(\pi^*H_r).
\end{equation}
We impose the sufficient condition
\begin{equation}\label{eq:resolution-nonzero-residues}
 \lambda_\nu\ne0,\qquad \nu=1,\ldots,N.
\end{equation}

\paragraph{The Poincar\'e metric.}
Fix a K\"ahler form $\omega_0$ on the smooth projective manifold $\widetilde X$.
Give each $\O(E_\nu)$ a smooth Hermitian metric and let $s_\nu$ be its canonical section.
Choose $C_\nu$ so that
\begin{equation*}
 \rho_\nu=C_\nu-\log\|s_\nu\|^2\ge1
 \quad\text{on }M,\qquad
 \Theta_\nu=\i\,\partial\bar\partial\rho_\nu.
\end{equation*}
Each $\Theta_\nu$ extends smoothly across $E$.
Compactness gives constants $K_\nu\ge0$ such that
$-K_\nu\omega_0\le\Theta_\nu\le K_\nu\omega_0$.
Positivity of $\omega_0$ on $\widetilde X$ makes these curvature bounds uniform and provides the ordinary-coordinate lower bound needed below.
The pullback $\pi^*\omega_{FS}$ is positive on $M$, but in the presence of exceptional fibres it degenerates at $E$ and does not supply this uniform background comparison.
For $A\ge1+\sum_\nu K_\nu$, define
\begin{equation}\label{eq:resolution-poincare-metric}
 \begin{aligned}
 \omega_P
 &=A\omega_0-\i\,\partial\bar\partial
                         \sum_{\nu=1}^N\log\rho_\nu\\
 &=\omega_B+\i\sum_{\nu=1}^N
       \frac{\partial\rho_\nu\wedge\bar\partial\rho_\nu}{\rho_\nu^2},
 \qquad
 \omega_B=A\omega_0-\sum_{\nu=1}^N\frac{\Theta_\nu}{\rho_\nu}.
 \end{aligned}
\end{equation}
Set $C_{B,0}:=A+\sum_\nu K_\nu$.
The background satisfies
\begin{equation*}
 \omega_0\le\omega_B\le C_{B,0}\omega_0.
\end{equation*}
Hence $\omega_P$ is positive and closed and defines a K\"ahler metric $g_P$.
This is the standard Poincar\'e type construction; see \parencite[Section~1.1]{Auvray2017Poincare}.

To prove completeness, set $U=\sum_\nu\log\rho_\nu$.
As a point approaches $E$, at least one $\rho_\nu$ tends to infinity, while all the others are bounded below by $1$.
Thus $U\to+\infty$ at the boundary.
Since $\widetilde X$ is compact, every curve escaping all compact subsets of $M$ approaches $E$.

For a piecewise $C^1$ curve $\gamma:[0,T)\to M$ escaping every compact subset, write $\dot\gamma=v+\bar v$ with $v\in T^{1,0}M$.
Then $\dot\gamma\log\rho_\nu=2\operatorname{Re}(\partial\rho_\nu(v)/\rho_\nu)$ and $g_P(\dot\gamma,\dot\gamma)=2g_P(v,v)$.
Together with \eqref{eq:resolution-poincare-metric}, these identities give
\begin{equation*}
 \left|\frac{\d}{\d t}\bigl(U(\gamma(t))\bigr)\right|^2
 \le N\sum_\nu|\dot\gamma\log\rho_\nu|^2
 \le2N g_P(\dot\gamma,\dot\gamma).
\end{equation*}
Consequently, for $0<t<T$, the arc length satisfies
\begin{equation*}
 \begin{aligned}
 \int_0^t\sqrt{g_P(\dot\gamma(s),\dot\gamma(s))}\,\d s
 &\ge\frac1{\sqrt{2N}}\int_0^t
       \left|\frac{\d}{\d s}\bigl(U(\gamma(s))\bigr)\right|\,\d s\\
 &\ge\frac1{\sqrt{2N}}
       |U(\gamma(t))-U(\gamma(0))|.
 \end{aligned}
\end{equation*}
The boundary behavior gives $U(\gamma(t))\to+\infty$ as $t\to T$, so every such curve has infinite arc length.
Hence $g_P$ is complete.

All estimate constants below are independent of the point and of the tangent vectors; they may depend on the fixed resolution, metrics, weights and coordinate neighborhood.
Set
\begin{equation*}
 R=\left(1+\sum_{\nu=1}^N\rho_\nu^2\right)^{1/2}.
\end{equation*}

Fix $p_0\in E$, and relabel the components so that precisely $E_1,\ldots,E_k$ contain $p_0$.
On a sufficiently small coordinate neighborhood, choose
\begin{equation*}
 E=\{z_1\cdots z_k=0\},\qquad
 e_i=z_i\partial_{z_i}\ (i\le k),\qquad
 e_j=\partial_{z_j}\ (j>k).
\end{equation*}
The dual holomorphic coframe on the complement is $\vartheta^i=\d z_i/z_i$ and $\vartheta^j=\d z_j$.
Work on a smaller neighborhood with compact closure in this coordinate chart.
Choose uniform constants $c_\rho,C_\rho>0$ so that, for $i\le k$,
\begin{equation*}
 t_i=-\log|z_i|^2\ge1,\qquad
 \rho_i=t_i+\psi_i,\qquad c_\rho t_i\le\rho_i\le C_\rho t_i,
\end{equation*}
where $\psi_i$ is smooth across $E$.
For $\nu>k$, set $\psi_\nu=\rho_\nu$, which is smooth since $s_\nu$ is nonvanishing on this neighborhood.
There are uniform constants $c_R,C_R>0$ such that
\begin{equation*}
 c_R\left(1+\sum_{i\le k}t_i^2\right)
 \le R^2\le
 C_R\left(1+\sum_{i\le k}t_i^2\right).
\end{equation*}

\begin{lemma}[Weighted comparison of the metric]\label{lem:poincare-weighted-metric}
Let $h$ be the product cusp metric
\begin{equation*}
 h=\sum_{i\le k}\frac{|\d z_i|^2}{|z_i|^2t_i^2}
                  +\sum_{j>k}|\d z_j|^2.
\end{equation*}
Define the positive scales by
\begin{equation*}
 \sigma_i=t_i^{-1}\quad(i\le k),\qquad
 \sigma_j=1\quad(j>k).
\end{equation*}
For $v=\sum_i v_i e_i$,
\begin{equation*}
 h(v,v)=\sum_{i\le k}\frac{|v_i|^2}{t_i^2}
          +\sum_{j>k}|v_j|^2
       =\sum_i\sigma_i^2|v_i|^2.
\end{equation*}
There exist constants $c_P,C_P>0$, independent of the point in the coordinate neighborhood intersected with $M$ and of $v$, such that
\begin{equation}\label{eq:resolution-model-comparison}
 c_P\sum_i\sigma_i^2|v_i|^2
 \le g_P(v,v)\le C_P\sum_i\sigma_i^2|v_i|^2.
\end{equation}
Equivalently, $c_Ph\le g_P\le C_Ph$.
\end{lemma}

\begin{proof}
Let $g_0,g_B$ be the Hermitian tensors associated with $\omega_0,\omega_B$, respectively, and put $\eta_\nu=\partial\rho_\nu/\rho_\nu$.
Equation~\eqref{eq:resolution-poincare-metric} gives
\begin{equation*}
 g_P(v,v)=g_B(v,v)+\sum_{\nu=1}^N|\eta_\nu(v)|^2,
 \qquad g_0\le g_B\le C_{B,0}g_0.
\end{equation*}
Put
\begin{equation*}
 B(v)=\sum_{i\le k}|z_i|^2|v_i|^2+\sum_{j>k}|v_j|^2.
\end{equation*}
This is the squared Euclidean norm of $v$ in the ordinary coordinates $z$.
The smooth positive definite metric $g_0$ is uniformly equivalent to this Euclidean metric on the compact closure.
Hence there are uniform constants $c_B,C_B>0$ such that
\begin{equation*}
 c_B\,B(v)\le g_B(v,v)\le C_B\,B(v).
\end{equation*}
For a uniform constant $C_{\partial\psi}>0$, the bounded ordinary coordinate derivatives of $\psi_\nu$ and Cauchy--Schwarz give $|\partial\psi_\nu(v)|\le C_{\partial\psi}B(v)^{1/2}$.
Since $\rho_\nu\ge1$, the constant $C_\psi:=NC_{\partial\psi}^2$ satisfies
\begin{equation*}
 \sum_{\nu=1}^N\frac{|\partial\psi_\nu(v)|^2}{\rho_\nu^2}
 \le C_\psi B(v).
\end{equation*}

The identities
\begin{equation*}
 \eta_i(v)=-\frac{v_i}{\rho_i}
              +\frac{\partial\psi_i(v)}{\rho_i}\quad(i\le k),
 \qquad
 \eta_\nu(v)=\frac{\partial\psi_\nu(v)}{\rho_\nu}\quad(\nu>k)
\end{equation*}
therefore imply
\begin{equation*}
 \begin{aligned}
 \sum_{\nu=1}^N|\eta_\nu(v)|^2
 &\le2\sum_{i\le k}\frac{|v_i|^2}{\rho_i^2}+2C_\psi B(v),\\
 \sum_{i\le k}\frac{|v_i|^2}{\rho_i^2}
 &\le2\sum_{\nu=1}^N|\eta_\nu(v)|^2+2C_\psi B(v).
 \end{aligned}
\end{equation*}
Together with the bounds for $g_B$, these yield
\begin{equation*}
 c_\eta\left(B(v)+\sum_{i\le k}\frac{|v_i|^2}{\rho_i^2}\right)
 \le g_P(v,v)\le
 C_\eta\left(B(v)+\sum_{i\le k}\frac{|v_i|^2}{\rho_i^2}\right),
\end{equation*}
where $c_\eta:=(2+(1+2C_\psi)/c_B)^{-1}$ and $C_\eta:=\max\{C_B+2C_\psi,2\}$.

The bounds for $\rho_i$ give
\begin{equation*}
 \frac{C_\rho^{-2}}{t_i^2}
 \le\frac1{\rho_i^2}\le\frac{c_\rho^{-2}}{t_i^2}.
\end{equation*}
Since $e^{-t}t^2$ is bounded for $t\ge1$, choose $K\ge1$ such that
\begin{equation*}
 |z_i|^2=e^{-t_i}\le\frac{K}{t_i^2},\qquad i\le k.
\end{equation*}
Thus $B(v)\le Kh(v,v)$.
For the upper bound,
\begin{equation*}
 B(v)+\sum_{i\le k}\frac{|v_i|^2}{\rho_i^2}
 \le(K+c_\rho^{-2})h(v,v)=C_h\,h(v,v),
\end{equation*}
where $C_h:=K+c_\rho^{-2}$.

For the lower bound, retaining the tangential part of $B(v)$ gives
\begin{equation*}
 \begin{aligned}
 B(v)+\sum_{i\le k}\frac{|v_i|^2}{\rho_i^2}
 &\ge\sum_{j>k}|v_j|^2
       +C_\rho^{-2}\sum_{i\le k}\frac{|v_i|^2}{t_i^2}\\
 &\ge c_h\,h(v,v),
 \end{aligned}
\end{equation*}
where $c_h:=\min\{1,C_\rho^{-2}\}$.
Combining the two comparisons proves the claim with $c_P=c_\eta c_h$ and $C_P=C_\eta C_h$.
\end{proof}

\begin{lemma}[Weighted connection estimates]\label{lem:poincare-connection}
Let $\nabla$ be the complexified Levi--Civita connection of $g_P$.
Let $c_P,C_P>0$ be comparison constants supplied by Lemma~\ref{lem:poincare-weighted-metric}.
In the holomorphic frame above, put
\begin{equation*}
 G_{i\bar j}=g_P(e_i,e_j),\qquad
 \nabla_{e_i}e_j=\Gamma^p_{ij}e_p.
\end{equation*}
There are uniform constants $C_{\partial G},C_\Gamma>0$ such that
\begin{equation*}
 \begin{aligned}
 |G^{\bar\ell p}|&\le\frac{c_P^{-1}}{\sigma_\ell\sigma_p},\\
 |e_iG_{j\bar\ell}|&\le C_{\partial G}\sigma_i\sigma_j\sigma_\ell,\\
 |\Gamma^p_{ij}|&\le C_\Gamma\frac{\sigma_i\sigma_j}{\sigma_p}.
 \end{aligned}
\end{equation*}
\end{lemma}

\begin{proof}
We first bound the inverse metric and record the induced tensor comparisons.
Inverting $c_Ph\le g_P\le C_Ph$ reverses the order, giving
\begin{equation*}
 C_P^{-1}h^{-1}\le g_P^{-1}\le c_P^{-1}h^{-1}.
\end{equation*}
Cauchy--Schwarz for the positive Hermitian matrix $G^{-1}$ gives
\begin{equation*}
 |G^{\bar\ell p}|
 \le\sqrt{G^{\bar\ell\ell}G^{\bar pp}}
 \le\frac{c_P^{-1}}{\sigma_\ell\sigma_p}.
\end{equation*}
Taking tensor products of the dual metrics gives, for every covariant tensor $T$ of order $r=1,2$,
\begin{equation}\label{eq:poincare-tensor-comparison}
 C_P^{-r/2}|T|_h\le|T|_{g_P}\le c_P^{-r/2}|T|_h.
\end{equation}

We next record the bounds needed to differentiate the metric coefficients.
Put $q_i=z_i$ for $i\le k$ and $q_i=1$ for $i>k$, so $e_i=q_i\partial_{z_i}$.
Set $K:=\max\{1,\sup_{t\ge1}e^{-t/2}t^2\}<\infty$; then, for $i\le k$,
\begin{equation*}
 |z_i|=e^{-t_i/2}\le\frac{K}{t_i^2}=K\sigma_i^2.
\end{equation*}
Since $0<\sigma_i\le1$ and $K\ge1$, this gives
\begin{equation*}
 |q_i|\le K\sigma_i,\qquad
 |e_iq_j|\le K\sigma_i\sigma_j.
\end{equation*}
Indeed, $e_iq_j$ vanishes unless $i=j\le k$, in which case it equals $z_i$.
For $\nu\le k$, the identities
\begin{equation*}
 e_jt_\nu=\bar e_jt_\nu=-\delta_{\nu j},\qquad
 e_i e_jt_\nu=e_i\bar e_jt_\nu=0
\end{equation*}
give $e_j\rho_\nu=-\delta_{\nu j}+q_j\partial_{z_j}\psi_\nu$.
For $\nu>k$, $e_j\rho_\nu=q_j\partial_{z_j}\psi_\nu$.
For every $\nu$, differentiating once more gives
\begin{equation*}
 \begin{aligned}
 e_i e_j\rho_\nu
 &=(e_iq_j)\partial_{z_j}\psi_\nu
       +q_iq_j\partial_{z_i}\partial_{z_j}\psi_\nu,\\
 e_i\bar e_j\rho_\nu
 &=q_i\bar q_j\partial_{z_i}\partial_{\bar z_j}\psi_\nu.
 \end{aligned}
\end{equation*}
The coordinate derivatives of $\psi_\nu$ are bounded on the compact closure.
Since $\rho_\nu^{-1}\le1$ for every $\nu$ and $\delta_{\nu j}/\rho_\nu\le c_\rho^{-1}\delta_{\nu j}\sigma_j$ for $\nu\le k$, the preceding identities and bounds give a uniform constant $C_{\partial\rho}>0$ such that
\begin{equation*}
 \begin{aligned}
 \frac{|e_j\rho_\nu|}{\rho_\nu}
 &\le C_{\partial\rho}\sigma_j,\\
 \frac{|e_i e_j\rho_\nu|+|e_i\bar e_j\rho_\nu|}{\rho_\nu}
 &\le C_{\partial\rho}\sigma_i\sigma_j.
 \end{aligned}
\end{equation*}

We now estimate the derivatives of the metric coefficients.
Let $g_0$ and $T_\nu$ be the Hermitian tensors associated with $\omega_0$ and $\Theta_\nu$, respectively.
Equation~\eqref{eq:resolution-poincare-metric} gives
\begin{equation*}
 G_{j\bar\ell}
 =Ag_0(e_j,e_\ell)-\sum_\nu\frac{T_\nu(e_j,e_\ell)}{\rho_\nu}
       +\sum_\nu\frac{(e_j\rho_\nu)(\bar e_\ell\rho_\nu)}{\rho_\nu^2}.
\end{equation*}
Differentiating yields
\begin{equation*}
 \begin{aligned}
 e_iG_{j\bar\ell}
 &=Ae_i[g_0(e_j,e_\ell)]\\
 &\quad+\sum_\nu\left[
       \frac{(e_i\rho_\nu)T_\nu(e_j,e_\ell)}{\rho_\nu^2}
       -\frac{e_i[T_\nu(e_j,e_\ell)]}{\rho_\nu}\right]\\
 &\quad+\sum_\nu\frac{(e_i e_j\rho_\nu)(\bar e_\ell\rho_\nu)
                    +(e_j\rho_\nu)(e_i\bar e_\ell\rho_\nu)}{\rho_\nu^2}\\
 &\quad-2\sum_\nu\frac{(e_i\rho_\nu)(e_j\rho_\nu)(\bar e_\ell\rho_\nu)}{\rho_\nu^3}.
 \end{aligned}
\end{equation*}
The ordinary-coordinate coefficients of $g_0,T_\nu$ and their first derivatives are bounded on the compact closure.
Writing $e_j=q_j\partial_{z_j}$ in the tensor evaluations and using $e_i\bar q_\ell=0$, their derivatives contain factors $(e_iq_j)\bar q_\ell$ or $q_iq_j\bar q_\ell$ multiplying these bounded coefficients.
Applying the preceding bounds to each term, with $\rho_\nu^{-1}\le1$ and $|\bar e_\ell\rho_\nu|=|e_\ell\rho_\nu|$, gives a uniform constant $C_{\partial G}>0$ such that
\begin{equation*}
 |e_iG_{j\bar\ell}|\le C_{\partial G}\sigma_i\sigma_j\sigma_\ell.
\end{equation*}

Finally, the K\"ahler connection formula gives
\begin{equation*}
 \Gamma^p_{ij}=(e_iG_{j\bar\ell})G^{\bar\ell p}.
\end{equation*}
The inverse and derivative bounds imply
\begin{equation*}
 |\Gamma^p_{ij}|
 \le\sum_\ell C_{\partial G}\sigma_i\sigma_j\sigma_\ell
                      \frac{c_P^{-1}}{\sigma_\ell\sigma_p}
 =nC_{\partial G}c_P^{-1}\frac{\sigma_i\sigma_j}{\sigma_p}.
\end{equation*}
Take $C_\Gamma:=nC_{\partial G}c_P^{-1}$.
\end{proof}

\begin{lemma}[Covariant derivative estimate]\label{lem:poincare-master-derivative}
There is a uniform constant $C_{\nabla\theta}>0$ such that
\begin{equation*}
 |\nabla\theta|_{g_P}\le C_{\nabla\theta}R
\end{equation*}
in a neighborhood of $E$.
\end{lemma}

\begin{proof}
We first estimate the coefficients of $\theta$ and their derivatives.
Using the logarithmic extension and the residues \eqref{eq:resolution-residues}, write $\theta=\sum_j\theta_j\vartheta^j$ with
\begin{equation*}
 \theta_j=
 \begin{cases}
 \lambda_j+z_j\beta_j,&j\le k,\\
 \beta_j,&j>k,
 \end{cases}
\end{equation*}
where $\beta=\sum_j\beta_j\d z_j$ extends holomorphically across $E$.
Choose $C_\beta>0$ bounding $|\beta_j|+|\partial_{z_i}\beta_j|$ for all $i,j$ on the compact closure.
Differentiating the coefficients gives
\begin{equation*}
 e_i\theta_j=(e_iq_j)\beta_j+q_je_i\beta_j.
\end{equation*}
The frame-weight bounds, with $K$ as in the preceding proof, yield
\begin{equation*}
 |e_i\theta_j|
 \le(K+K^2)C_\beta\sigma_i\sigma_j
 =C_{\partial\theta}\sigma_i\sigma_j,
 \qquad C_{\partial\theta}:=(K+K^2)C_\beta.
\end{equation*}
Set $C_{\theta,h}:=c_R^{-1/2}[(\sum_{i\le k}|\lambda_i|^2)^{1/2}+\sqrt nKC_\beta]$.
The diagonal expression for $h$ and the lower bound for $R^2$ give
\begin{equation*}
 \begin{aligned}
 |\theta|_h
 &=\left(\sum_j\frac{|\theta_j|^2}{\sigma_j^2}\right)^{1/2}\\
 &\le\left(\sum_{i\le k}|\lambda_i|^2t_i^2\right)^{1/2}
        +\sqrt nKC_\beta\\
 &\le C_{\theta,h}R.
 \end{aligned}
\end{equation*}

We next estimate the covariant derivative in the logarithmic frame.
Its components satisfy
\begin{equation*}
 (\nabla_{e_i}\theta)(e_j)
 =e_i\theta_j-\sum_p\Gamma^p_{ij}\theta_p,\qquad
 \nabla_{\bar e_i}\theta=0,
\end{equation*}
where the second identity uses holomorphicity and the K\"ahler property.
Since $|\theta_p|\le\sigma_p|\theta|_h\le C_{\theta,h}\sigma_pR$, the preceding lemma gives
\begin{equation*}
 \begin{aligned}
 |(\nabla_{e_i}\theta)(e_j)|
 &\le C_{\partial\theta}\sigma_i\sigma_j
       +\sum_p C_\Gamma\frac{\sigma_i\sigma_j}{\sigma_p}
                     C_{\theta,h}\sigma_pR\\
 &\le(C_{\partial\theta}+nC_\Gamma C_{\theta,h})
                     \sigma_i\sigma_jR,
 \end{aligned}
\end{equation*}
since $R\ge1$.

Finally, \eqref{eq:poincare-tensor-comparison} with $r=2$ gives
\begin{equation*}
 \begin{aligned}
 |\nabla\theta|_{g_P}
 &\le c_P^{-1}|\nabla\theta|_h\\
 &=c_P^{-1}\left(\sum_{i,j}
        \frac{|(\nabla_{e_i}\theta)(e_j)|^2}
             {\sigma_i^2\sigma_j^2}\right)^{1/2}\\
 &\le\frac{n}{c_P}(C_{\partial\theta}+nC_\Gamma C_{\theta,h})R.
 \end{aligned}
\end{equation*}
Taking the maximum of these coefficients over a finite boundary cover proves the claim.
\end{proof}

\begin{lemma}[Norm estimates]\label{lem:poincare-master-norm}
If every $\lambda_\nu\ne0$, there exist $K\Subset M$ and $c_{\theta},C_{\theta}>0$ such that
\begin{equation*}
 c_{\theta} R^2\le|\theta|_{g_P}^2\le C_{\theta} R^2
 \qquad\text{on }M\setminus K.
\end{equation*}
\end{lemma}

\begin{proof}
Let $\lambda_*:=\min_\nu|\lambda_\nu|>0$ and shrink the boundary neighborhoods so that $|z_i\beta_i|\le\lambda_*/2$ for $i\le k$.
Then
\begin{equation*}
 \begin{aligned}
 |\theta|_h^2
 &=\sum_{i\le k}t_i^2|\lambda_i+z_i\beta_i|^2
                         +\sum_{j>k}|\beta_j|^2\\
 &\ge\frac{\lambda_*^2}{4}\sum_{i\le k}t_i^2.
 \end{aligned}
\end{equation*}
Since $k\ge1$ and $t_i\ge1$, the upper bound for $R^2$ gives $R^2\le2C_R\sum_{i\le k}t_i^2$.
The bound $|\theta|_h\le C_{\theta,h}R$ from the preceding proof and \eqref{eq:poincare-tensor-comparison} with $r=1$ give
\begin{equation*}
 \begin{aligned}
 |\theta|_{g_P}^2
 &\ge C_P^{-1}|\theta|_h^2
 \ge\frac{\lambda_*^2}{8C_PC_R}R^2,\\
 |\theta|_{g_P}^2
 &\le c_P^{-1}|\theta|_h^2
 \le c_P^{-1}C_{\theta,h}^2R^2.
 \end{aligned}
\end{equation*}
Choose a finite collection of smaller boundary coordinate neighborhoods.
Take $c_{\theta}$ to be the minimum of their coefficients $\lambda_*^2/(8C_PC_R)$ and $C_{\theta}$ to be the maximum of their coefficients $c_P^{-1}C_{\theta,h}^2$.
The required thresholds can also be made uniform.
The complement of their union is compact in $M$, proving the claim.
\end{proof}

\begin{proposition}\label{prop:well-tamed-cylindrical}
Assume that the boundary residues satisfy \eqref{eq:resolution-nonzero-residues}.
For the complete metric $g_P$ of \eqref{eq:resolution-poincare-metric}, there are constants $c_{\theta},C_{\theta},C_{\nabla\theta}>0$ and a compact set $K\Subset M$ such that
\begin{equation}\label{eq:resolution-poincare-estimates}
 c_{\theta} R^2\le|\theta|_{g_P}^2\le C_{\theta} R^2,\qquad
 |\nabla^{g_P}\theta|_{g_P}\le C_{\nabla\theta}R
 \quad\text{on }M\setminus K.
\end{equation}
In particular, $(M,g_P,\theta)$ is well-tamed in the sense of Definition~\ref{def:well-tamed}.
\end{proposition}

\begin{proof}
The form $\theta$ is closed and holomorphic on $M$ and is nontrivial because its residues are nonzero.
Completeness was proved above.
The two preceding lemmas give \eqref{eq:resolution-poincare-estimates} outside a compact subset.
Since $R\to\infty$ on every sequence escaping compact subsets, enlarge $K$ so that $c_{\theta} R^2>1$ and $R>2C_{\nabla\theta}/c_{\theta}$ on $M\setminus K$.
Set $C_1=C_2=1$.
Then $|\theta|_{g_P}>C_1$ and
\begin{equation*}
 \begin{aligned}
 |\theta|_{g_P}^2-C_2|\nabla^{g_P}\theta|_{g_P}
 &\ge c_{\theta} R^2-C_{\nabla\theta}R\\
 &\ge\frac{c_{\theta}}{2}R^2>0,
 \end{aligned}
\end{equation*}
which proves well-tameness.
\end{proof}

On the maximal wonderful model, the boundary components are indexed by the nonempty proper intersections $X$ of the arrangement hyperplanes.
At a general point of $X$, each defining linear section vanishes to order one exactly when $X\subset H_r$.
The corresponding boundary valuation consequently gives
\begin{equation*}
 \lambda_X=\sum_{H_r\supset X}a_r=a_X.
\end{equation*}
Thus the nonresonance condition \eqref{eq:nonresonance-weights} implies \eqref{eq:resolution-nonzero-residues} on this model.
In this case, the counting result in Section~\ref{sec:YY} gives the total Milnor number by \eqref{eq:counting-lattice}.
It equals the number of distinct critical points when all critical points are nondegenerate.

\section{The Yang--Yang functional and the Bethe equations}\label{sec:YY-functional}

We apply Sections~\ref{sec:main_theorem} and \ref{sec:YY} to the logarithmic differential of the Yang--Yang functional of the Gaudin model \parencite{NEKRASOV201169,Liu_2021,HuLiu2015}.
For fixed marked points, its critical-point equations are the Bethe equations.
After a log resolution of the projective collision arrangement, nonzero residues on every boundary component allow us to apply the Poincar\'e metric construction of Section~\ref{sec5.4}.

Let $\mathfrak g$ be a semisimple Lie algebra with simple roots $\alpha_1,\ldots,\alpha_r$.
Let $\mathfrak h^*$ be the dual of a Cartan subalgebra with invariant bilinear form $(\, ,\, )$.
Fix $m,N\ge1$, dominant integral weights $\Lambda_1,\ldots,\Lambda_N$, distinct marked points $z_1,\ldots,z_N\in\C$, and a color map $c:\{1,\ldots,m\}\to\{1,\ldots,r\}$.
The ordered Bethe variables $t=(t_1,\ldots,t_m)$ range over
\begin{equation*}
 M_z=\{t\in\C^m:t_j\ne z_a\text{ for all }j,a,\quad
                         t_i\ne t_j\text{ for }i\ne j\}.
\end{equation*}
The Yang--Yang functional is the multivalued function
\begin{equation}\label{eq:YY-functional}
 \begin{aligned}
 \Phi(t;z)
 &=\sum_{1\le a<b\le N}(\Lambda_a,\Lambda_b)\log(z_a-z_b)\\
 &\quad-\sum_{j=1}^m\sum_{a=1}^N
             (\alpha_{c(j)},\Lambda_a)\log(t_j-z_a)\\
 &\quad+\sum_{1\le i<j\le m}
             (\alpha_{c(i)},\alpha_{c(j)})\log(t_i-t_j).
 \end{aligned}
\end{equation}
Throughout this section, $z$ is fixed and $\theta_z:=\d_t\Phi(t;z)$.
The first sum in \eqref{eq:YY-functional} is independent of $t$ and can be omitted when computing $\theta_z$.
The Bethe equations $\theta_z=0$ are
\begin{equation}\label{eq:Bethe-equations}
 -\sum_{a=1}^N\frac{(\alpha_{c(j)},\Lambda_a)}{t_j-z_a}
 +\sum_{i\ne j}\frac{(\alpha_{c(j)},\alpha_{c(i)})}{t_j-t_i}
 =0,\qquad j=1,\ldots,m.
\end{equation}
For $\mathfrak g=\mathfrak{sl}_2$, normalize $(\alpha,\alpha)=2$ and write $m_a=(\alpha,\Lambda_a)$.
Up to a term depending only on $z$, the functional becomes
\begin{equation*}
 \Phi(t;z)=-\sum_{j=1}^m\sum_{a=1}^N m_a\log(t_j-z_a)
                 +2\sum_{1\le i<j\le m}\log(t_i-t_j),
\end{equation*}
and the Bethe equations are
\begin{equation}\label{eq:sl2-Bethe}
 -\sum_{a=1}^N\frac{m_a}{t_j-z_a}
 +2\sum_{i\ne j}\frac1{t_j-t_i}=0,
 \qquad j=1,\ldots,m.
\end{equation}

In homogeneous coordinates $[T_0:\cdots:T_m]$ on $\P^m$, let
\begin{equation*}
 H_{j,a}=\{T_j-z_aT_0=0\},\qquad
 H_{ij}=\{T_i-T_j=0\},\qquad H_\infty=\{T_0=0\}.
\end{equation*}
Write $\mathcal A_z=\{H_\infty\}\cup\{H_{j,a}\}_{j,a}\cup\{H_{ij}\}_{i<j}$ and
\begin{equation*}
 D_z=\bigcup_{H\in\mathcal A_z}H,\qquad M_z=\P^m\setminus D_z.
\end{equation*}
For $m\ge2$, the collision arrangement is not simple normal crossing: in $\C^2$, the three lines $t_1=z_1$, $t_2=z_1$, and $t_1=t_2$ meet at the point $t_1=t_2=z_1$, which has codimension two.
Assign the finite hyperplanes the weights
\begin{equation*}
 b_{j,a}=-(\alpha_{c(j)},\Lambda_a),\qquad
 b_{ij}=(\alpha_{c(i)},\alpha_{c(j)}),
\end{equation*}
and assign $H_\infty$ the weight
\begin{equation*}
 b_\infty:=-\sum_{j,a}b_{j,a}-\sum_{i<j}b_{ij}.
\end{equation*}
For $H\in\mathcal A_z$, write $b_H$ for its assigned weight.
These are the residues of $\theta_z$ along the corresponding projective hyperplanes.

Choose the maximal De Concini--Procesi wonderful model
\begin{equation*}
 \pi:\widetilde X_z\longrightarrow\P^m,\qquad
 E=(\pi^{-1}D_z)_{\mathrm{red}}=\bigcup_\nu E_\nu.
\end{equation*}
Then $E$ is simple normal crossing and $\pi$ is an isomorphism over $M_z$.
Identify $M_z$ with $\widetilde X_z\setminus E$ and retain the notation $\theta_z$ for the pullback.
The form $\theta_z$ extends logarithmically across $E$, with residues
\begin{equation}\label{eq:YY-boundary-residues}
 \lambda_\nu=\sum_{H\in\mathcal A_z}
              b_H\operatorname{ord}_{E_\nu}(\pi^*H).
\end{equation}
We now describe the boundary components and compute their residues to state a sufficient condition for applying Section~\ref{sec5.4}.
On this model, each nonempty proper intersection $X$ of the arrangement hyperplanes determines one irreducible boundary component $E_X$.
For a hyperplane $X$, $E_X$ is its strict transform.
Intersections of codimension at least two give additional boundary components, including intersections where the hyperplanes already meet normally.

Since the finite hyperplanes have equations $t_j=z_a$ or $t_i=t_j$, their intersections are specified by groups of variables with prescribed or common values.
The resulting boundary components are described by groups approaching marked points, common unmarked finite values, or infinity, with the equalities at infinity understood after rescaling.

We begin with boundary components for variables approaching a marked point.
For a fixed marked point $z_a$ and a nonempty index set $S\subset\{1,\ldots,m\}$, the corresponding intersection is
\begin{equation*}
 X_{S,a}=\{T_j=z_aT_0:j\in S\}.
\end{equation*}
The residue is
\begin{equation*}
 \lambda_{S,a}=\sum_{j\in S}b_{j,a}
       +\sum_{\substack{i<j\\i,j\in S}}b_{ij},
 \qquad |S|\ge1.
\end{equation*}

We next consider boundary components for variables approaching a common unmarked finite value.
Let $S\subset\{1,\ldots,m\}$ index these variables, with $|S|\ge2$; the corresponding intersection is
\begin{equation*}
 X_S=\{T_i=T_j:i,j\in S\}.
\end{equation*}
The residue is
\begin{equation*}
 \lambda_S=\sum_{\substack{i<j\\i,j\in S}}b_{ij},\qquad |S|\ge2.
\end{equation*}

We finally consider boundary components for variables tending to infinity with distinct nonzero limits after rescaling.
For a nonempty index set $S\subset\{1,\ldots,m\}$, the corresponding intersection is
\begin{equation*}
 X_{S,\infty}=\{T_0=0,\ T_j=0\text{ for }j\notin S\}.
\end{equation*}
The hyperplanes containing $X_{S,\infty}$ are $H_\infty$, the $H_{j,a}$ with $j\notin S$, and the $H_{ij}$ with $i,j\notin S$.
Thus \eqref{eq:YY-boundary-residues} gives
\begin{equation}\label{eq:YY-infinity-residue}
 \begin{aligned}
 \lambda_{S,\infty}
 &=b_\infty+\sum_{j\notin S}\sum_a b_{j,a}
              +\sum_{\substack{i<j\\i,j\notin S}}b_{ij}\\
 &=-\sum_{j\in S}\sum_a b_{j,a}
              -\sum_{\substack{i<j\\\{i,j\}\cap S\ne\varnothing}}b_{ij}.
 \end{aligned}
\end{equation}
For $S=\{1,\ldots,m\}$, this is $b_\infty$.

\begin{definition}\label{def:YY-nonresonant}
The weights of the Yang--Yang functional are called nonresonant on the chosen model if $\lambda_\nu\ne0$ for every irreducible boundary component $E_\nu$.
\end{definition}

On the maximal wonderful model, this is the nonresonance condition \eqref{eq:nonresonance-weights} for the projective arrangement $\mathcal A_z$.
It is an additional restriction on the Lie-theoretic data.
For example, in the $\mathfrak{sl}_2$ case, $b_{j,a}=-m_a$ and $b_{ij}=2$, so a marked collision with $s=|S|$ has $\lambda_{S,a}=s(s-1-m_a)$.
Thus nonresonance fails if $m\ge m_a+1$ for some marked point $z_a$.

\begin{proposition}\label{prop:YY-finiteness}
Assume the nonresonance condition of Definition~\ref{def:YY-nonresonant}.
Let $g_{P,z}$ be the complete Poincar\'e metric constructed in Section~\ref{sec5.4}.
Then $(M_z,g_{P,z},\theta_z)$ is a well-tamed Novikov type LG model.
There are noncanonical vector-space isomorphisms
\begin{equation*}
 \mathcal H^q_{\theta_z}\cong
 \begin{cases}
  \displaystyle\bigoplus_{p\in\mathrm{Crit}(\theta_z)}
                      \mathcal R_{\theta_z,p},&q=m,\\
  0,&q\ne m,
 \end{cases}
\end{equation*}
and
\begin{equation*}
 \begin{aligned}
 \dim_\C\mathcal H^m_{\theta_z}
 &=\sum_{p\in\mathrm{Crit}(\theta_z)}\mu_p\\
 &=(-1)^m\sum_{X\in L(\mathcal A_z)}\mu(X)(\dim X+1).
 \end{aligned}
\end{equation*}
\end{proposition}

\begin{proof}
Proposition~\ref{prop:well-tamed-cylindrical} gives well-tameness for $g_{P,z}$.
The complement $M_z$ is a principal affine open subset of $\C^m$, hence is Stein, so Theorem~\ref{main_theorem_1} applies.
The affine volume form trivializes the canonical bundle of $M_z$, giving the stated vector-space identification with the local Milnor rings.
The dimension formula follows from \eqref{eq:counting-lattice} for $\mathcal A_z$.
\end{proof}

If all critical points are nondegenerate, each $\mu_p=1$, and the dimension formula counts distinct solutions in the ordered configuration space $M_z$.
The Bethe equations are invariant under permutations of variables of the same color \parencite[Section~2.1]{Liu_2021}.
For $n_\ell:=\#\{j:c(j)=\ell\}$, these permutations act freely on $M_z$, so the number of such solutions modulo same-color permutations is the ordered count divided by $\prod_{\ell=1}^r n_\ell!$.

\printbibliography

\end{document}